\documentclass[11pt]{article}
\usepackage{amsmath}
\usepackage{amsthm}
\usepackage{amsfonts}
\usepackage{amssymb}
\usepackage{latexsym}
\usepackage{diagbox}
\usepackage{mathtools}
\usepackage{xcolor}

\usepackage{array}

\usepackage{subcaption}
\usepackage{booktabs}
\usepackage{tkz-graph}
\usepackage{graphics,graphicx}
\usepackage{multirow}

\usepackage{tikz}
\usepackage{tikz-cd}

\usetikzlibrary{matrix,arrows}
\usetikzlibrary{fit}
\usetikzlibrary{patterns}
\usetikzlibrary{chains,fit,shapes}
\usetikzlibrary{positioning,calc}
\usetikzlibrary{graphs}
\usetikzlibrary{arrows,decorations.markings}
\usepackage[colorlinks,
linkcolor=blue,
anchorcolor=blue,
citecolor=blue
]{hyperref}

\newtheorem{thm}{Theorem}

\newtheorem{cor}[thm]{Corollary}

\newtheorem{rem}[thm]{Remark}

\DeclareMathOperator{\asc}{asc}
\DeclareMathOperator{\des}{des}
\DeclareMathOperator{\rmax}{rmax}
\DeclareMathOperator{\lmax}{lmax}
\DeclareMathOperator{\rmin}{rmin}
\DeclareMathOperator{\lmin}{lmin}

\begin{document}
\begin{center}
{\large \bf  Distribution of statistics on separable permutations restricted by a flat POP}
\end{center}
\begin{center}
Alice L.L. Gao$^{1}$, Sergey Kitaev$^{2}$, Ya-Xing Li$^{3}$ and Xuan Ruan$^{4}$\\[6pt]

$^{1,3,4}$Shenzhen Research Institute of Northwestern Polytechnical University,\\
Sanhang Science $\&$ Technology Building, No. 45th, Gaoxin South 9th Road, Nanshan District,   Shenzhen City,  518057,  P.R. China

$^{1,3,4}$School of Mathematics and Statistics,\\
 Northwestern Polytechnical University, Xi'an, Shaanxi 710072, P.R. China

 $^{2}$Department of Mathematics, \\
Department of Mathematics and Statistics, University of Strathclyde, 26 Richmond Street, Glasgow G1 1XH, United Kingdom \\[6pt]

 Email:$^{1}${\tt llgao@nwpu.edu.cn},
        $^{2}${\tt sergey.kitaev@strath.ac.uk},
       $^{3}${\tt liya@mail.nwpu.edu.cn},
       $^{4}${\tt ruanxuan1211@163.com}.
\end{center}

\noindent\textbf{Abstract.}
Finding distributions of statistics in pattern-avoiding permutations has attracted  significant attention in the literature. In particular, Chen, Kitaev, and Zhang derived functional equations for the joint distributions of any subset of classical minima and maxima statistics, as well as for the joint distributions of ascents and descents in separable permutations. Meanwhile, partially ordered patterns (POPs) have also been extensively studied. Notably, so-called flat POPs played a key role, via the notion of shape-Wilf-equivalence, in proving a conjecture on pattern-avoiding permutations.

In this paper, we study flat POP-avoiding separable permutations, where the maximum element in a flat POP receives the largest label. Avoiding such a POP imposes restrictions on the position of the maximum element in a separable permutation, forcing it to be positioned to the left. We establish a system of functional equations describing the joint distribution of six classical statistics in the most general case, extending the work of Chen, Kitaev, and Zhang.

As a specialization, when the POP has length 3, we recover a joint distribution result of Han and Kitaev on permutations avoiding classical patterns of length 3. As another specialization, for the flat POP of length 4, we derive an explicit rational generating function that captures the distribution of six statistics, with a numerator containing 100 monomials and a denominator containing 19 monomials.  \\

\noindent {\bf AMS Classification 2020:} 05A15
	
	
\section{Introduction}\label{intro-sec}

A permutation of length $n$ is a rearrangement of the set $[n]:=\{1,2,\ldots,n\}$. Denote by  $S_n$  the set of permutations of $[n]$ and let $\varepsilon$ be the empty permutation. A {\em pattern} is a permutation. A permutation $\pi_1\pi_2\cdots\pi_n\in S_n$ avoids a pattern $p=p_1p_2\cdots p_k\in S_k$  if there is no subsequence $\pi_{i_1}\pi_{i_2}\cdots\pi_{i_k}$ such that $\pi_{i_j}<\pi_{i_m}$ if and only if $p_j<p_m$. For example, the permutation $32154$ avoids the pattern $231$. Let $S_n(p)$ denote the set of $p$-avoiding permutations of length $n$. For a permutation $\pi=\pi_1\pi_2\cdots\pi_n$, its {\em reverse} is the permutation $r(\pi)=\pi_n\pi_{n-1}\cdots \pi_1$ and its {\em complement} is the permutation $c(\pi)=c(\pi_1)c(\pi_2)\cdots c(\pi_n)$ where $c(x)=n+1-x$. Also, the length of a permutation $\pi$, denoted by $|\pi|$, is the number of elements in $\pi$. For example, for $\pi=423165$, $|\pi|=6$, $r(\pi)=561324$, and $c(\pi)=354612$.
	
\begin{figure}[htbp]
		\centering
		\begin{tikzpicture}[scale=0.8]
			
			\draw [line width=1](0,0)--(1.5,1.5);
			\draw [line width=1](1,0)--(1.5,1.5);
			\draw [line width=1](3,0)--(1.5,1.5);

			\draw (0,0) node [scale=0.4, circle, draw,fill=black]{};
			\draw (1,0) node [scale=0.4, circle, draw,fill=black]{};
			\draw (3,0) node [scale=0.4, circle, draw,fill=black]{};
			\draw (1.5,1.5) node [scale=0.4, circle, draw,fill=black]{};
			
			\draw (1.75,0) node [scale=0.15, circle, draw,fill=black]{};
			\draw (2,0) node [scale=0.15, circle, draw,fill=black]{};
			\draw (2.25,0) node [scale=0.15, circle, draw,fill=black]{};
			
			
			\node [below] at (0,-0.1){$1$};
			\node [below] at (1,-0.1){$2$};
			\node [below] at (3,-0.1){${k-1}$};
			\node [above] at (1.5,1.5){${k}$};
		\end{tikzpicture}
		\vspace{-0.3cm}
		\caption{The form of POPs in this paper.}
		\label{pic-B5}
	\end{figure}
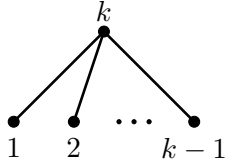

	A \emph{partially ordered pattern} ({POP}) $p$ of length $k$ is defined by a $k$-element partially ordered set (poset) $P$ labeled by the elements in $\{1,\ldots,k\}$. An occurrence of such a POP $p$ in a permutation $\pi=\pi_1\cdots\pi_n$ is a subsequence $\pi_{i_1}\cdots\pi_{i_k}$, where $1\le i_1<\cdots< i_k\le n$,  such that $\pi_{i_j}<\pi_{i_m}$ if and only if $j<m$ in $P$. Thus, a classical pattern of length $k$ corresponds to a $k$-element chain. For example, the POP $p=$ \hspace{-3.5mm}
	\begin{minipage}[c]{3.5em}
		\scalebox{1}{
			\begin{tikzpicture}[scale=0.5]
				
				\draw [line width=1](0,-0.5)--(0,0.5);
				
				\foreach \x/\y in {0/-0.5,1/-0.5,0/0.5}
				\draw (\x,\y) node [scale=0.4, circle, draw, fill=black]{};
				
				\node [left] at (0,-0.6){${\small 3}$};
				\node [right] at (1,-0.6){${\small 2}$};
				\node [left] at (0,0.6){${\small 1}$};
				
			\end{tikzpicture}
		}
	\end{minipage}
	occurs six times in the permutation $41523$, namely, as the subsequences $412$, $413$, $452$, $453$, $423$, and $523$. Clearly, avoiding $p$ is the same as avoiding the patterns $312$, $321$, and $231$ at the same time. In this paper, for $k \geq 1$, let $P_k$ denote a POP of the form shown in Figure~\ref{pic-B5}, and let $S_n(2413,3142,P_k)$ be the set of permutations in $S_n$ that avoid the patterns 2413, 3142, and $P_k$. By definition, $P_1 = 1$, the permutation of length one.

\vspace{-0.5cm}

Suppose $\pi=\pi_1\pi_2\cdots\pi_m \in S_m$ and $\sigma=\sigma_1\sigma_2\cdots\sigma_n\in S_n$. We define the {\em direct sum}\index{direct sum $\oplus$} (or simply, {\em sum}) $\oplus$, and the {\em skew sum}\index{skew sum $\ominus$} $\ominus$ by building the permutations $\pi\oplus\sigma$ and $\pi\ominus\sigma$ as follows:
\begin{eqnarray*}
		(\pi\oplus\sigma)_i&=& \left\{
		\begin{array}{ll}
			\pi_i  &  \mbox{ if }1\leq i\leq m,\\
			\sigma_{i-m}+m & \mbox{ if }m+1\leq i\leq m+n,
		\end{array}\right.\nonumber \\
		(\pi\ominus\sigma)_i &=& \left\{
		\begin{array}{ll}
			\pi_i+n  &  \mbox{ if }1\leq i\leq m,\\
			\sigma_{i-m} & \mbox{ if }m+1\leq i\leq m+n.
		\end{array}\right.\nonumber
	\end{eqnarray*}
	For example, $14325\oplus 4231=143259786$ and $14325\ominus 4231= 587694231$.
	The {\em separable permutations} are those which can be built from the permutation $1$ by repeatedly applying the $\oplus$ and $\ominus$ operations.  By definition, we assume that the empty permutation $\varepsilon$ is a separable permutation, though some papers make the opposite assumption. Bose, Buss and Lubiw~\cite{BBL98} introduced the notion of separable permutation in 1998, and it is well-known \cite{Kitaev2011Patterns} that the set of all separable permutations of length $n\geq 1$ is precisely $S_n(2413,3142)$.  Separable permutations appear in the literature in various contexts (see, for example, \cite{AAV2011,AHP2015,AJ2016,FuLinZeng,GaoLiu,NRV,Stankova}) and they are shown in \cite{West1995} to be counted by the $(n-1)$-th {\em Schr\"oder number}. A particular result of Chen, Kitaev and Zhang \cite{DAM2024CZS} is that the generating function (g.f.)
$$\sum_{n \geq 1} x^n \sum_{\pi\in S_n(2413, 3142)} u^{\text{lmax}(\pi)}$$
is given by
\begin{footnotesize}
\begin{equation}\label{F(x,u)}
\frac{4\sqrt{\left(-\frac{1}{4} \sqrt{x^2-6 x+1}-xu+\frac{x}{4}+\frac{5}{4}\right)^2+\sqrt{x^2-6 x+1}-x-1}-\sqrt{x^2-6 x+1}+4x u+x - 3}{2 \left(\sqrt{x^2-6 x+1}-x-1\right)}.
\end{equation}
\end{footnotesize}

\begin{figure}[htbp]
		\begin{center}
			\begin{tikzpicture}[line width=0.5pt,scale=0.24]
				\coordinate (O) at (0,0);
				
				\path (25,1)  node {$n$};
				\draw [dashed] (O)--++(50,0);
				\fill[black!100] (O)++(25,0) circle(1.5ex);
				
				\draw (20,-1) rectangle (24,-3);
				\path (22,-2)  node {$L_m$};
				\path (13.5,-2)  node {possibly empty$\boldsymbol{\rightarrow}$};
				
				\draw [dashed] (0,-4)--++(50,0);
				\draw [dashed] (0,-8)--++(50,0);
				\draw [dashed] (0,-17)--++(50,0);
				\draw [dashed] (0,-21)--++(50,0);
				\draw [dashed] (0,-25)--++(50,0);
				
				\draw [dashed] (25,0)--++(0,-29);
				\draw [dashed] (31,0)--++(0,-29);
				\draw [dashed] (37,0)--++(0,-29);
				\draw [dashed] (43,0)--++(0,-29);
				
				\draw [dashed] (19,0)--++(0,-29);
				\draw [dashed] (13,0)--++(0,-29);
				\draw [dashed] (7,0)--++(0,-29);
				
				\draw (26,-5) rectangle (30,-7);
				\path (28,-6)  node {$R_m$};
				
				\fill[black!100] (O)++(17,-9) circle(1.0ex);
				\fill[black!100] (O)++(16,-10) circle(1.0ex);
				\fill[black!100] (O)++(15,-11) circle(1.0ex);
				\draw (8,-14) rectangle (12,-16);
				\path (10,-15)  node {$L_2$};
				
				\draw (2,-22) rectangle (6,-24);
				\path (4,-23)  node {$L_1$};

				\fill[black!100] (O)++(32,-11) circle(1.0ex);
				\fill[black!100] (O)++(33,-12) circle(1.0ex);
				\fill[black!100] (O)++(34,-13) circle(1.0ex);
				\draw (38,-18) rectangle (42,-20);
				\path (40,-19)  node {$R_2$};
				
				\draw (45,-26) rectangle (49,-28);
				\path (47,-27)  node {$R_1$};
				\path (38.5,-27)  node {possibly empty$\boldsymbol{\rightarrow}$};
				
			\end{tikzpicture}
			\caption{A schematic representation of the permutation diagrams corresponding to separable permutations, known as \emph{Stankova's decomposition of separable permutations}~\cite{Stankova}. Each $L_i$ and $R_j$ is a separable permutation.}\label{sepStructure}
			\label{fig:stankova}
		\end{center}
	\end{figure}
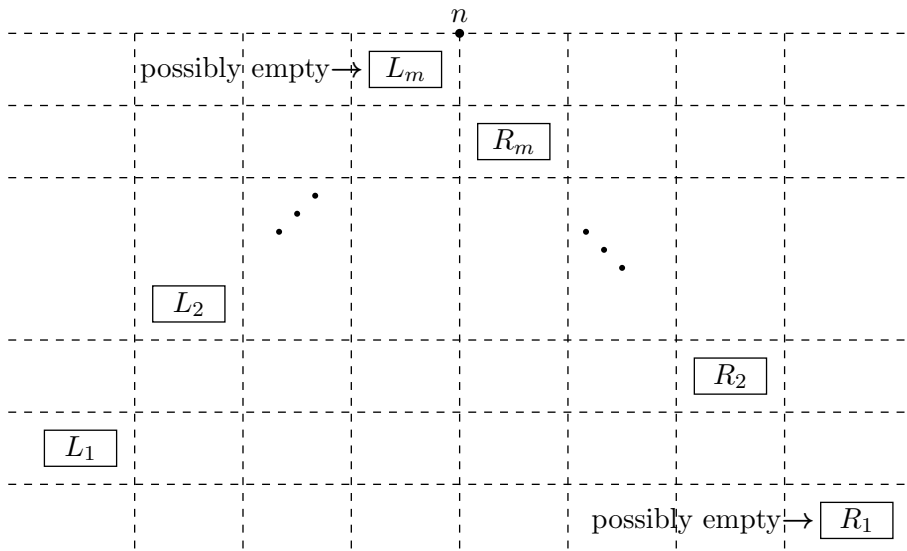

Clearly, the reverse, or complement, or (usual group theoretic) inverse, or any composition of these operations applied to a separable permutation gives a separable permutation. According to the seminal work of Stankova~\cite{Stankova}, any $\pi\in S_n(2413,3142)$ has the following structure
 (see also Figure~\ref{sepStructure} for a schematic representation, where a permutation is viewed as a diagram with a dot in position $(i,\pi_i)$ for each element $\pi_i$ of $\pi$):
\begin{align}
\pi=L_1L_2\cdots
		L_mnR_mR_{m-1}\cdots R_1\notag
\end{align}
where
\begin{itemize}
	\item[(1)] for $1\leq i\leq m$, $L_i$ and $R_i$ are non-empty ($\neq \varepsilon$), with possible exception for
	$L_m$ and $R_1$, separable permutations which are intervals in $\pi$ (that is, consist of all elements in $\{a,a+1,\ldots,b\}$ for some $a$ and $b$);
   \item[(2)]  $R_1<L_1<R_2<L_2<\cdots <R_m<L_m$, where $A<B$, for two
		permutations $A$ and $B$, means that each element of $A$ is less than
		every element of $B$. In particular, $R_1$, if it is not empty,
		contains 1.
	\end{itemize}
	
\noindent
For example, if $\pi=2165743$ then $L_1=21$, $L_2=65$, $R_1=\varepsilon$ and $R_2=43$.

Stankova's decomposition of separable permutations has proven to be highly effective in studying classical statistics on separable permutations. In particular, Lin~\cite{Lin2017} used the decomposition to compute the joint distribution of descents and double descents, providing an algebraic proof of a $\gamma$-positivity expansion originally due to Fu-Lin-Zeng~\cite{FuLinZeng}. A similar approach was employed by Wang~\cite{Wang} to establish a further generalization. Furthermore, the Comtet statistics \emph{comp} and \emph{iar}, denoting the number of components and the length of the initial ascending run, respectively, were shown to be symmetric over separable permutations~\cite{FuLinWang2020,FuLinWang2020b}. Most recently, Chen, Kitaev, and Zhang~\cite{DAM2024CZS} derived generating functions for the distributions of Eulerian and Stirling statistics on separable permutations. These results underscore the continuing significance of structural decompositions in the enumerative study of pattern-avoiding permutations.
	
In this paper, we study separable permutations that avoid the pattern $P_k$, i.e., the set $S_n(2413, 3142,P_k)$.
We focus on the following classical permutation statistics. For \( 1 \leq i \leq n - 1 \), if \( \pi_i < \pi_{i+1} \) (resp., \( \pi_i > \pi_{i+1} \)), then \( i \) is an {\em ascent} (resp., {\em descent}) in \( \pi \in S_n \). The number of ascents (resp., descents) in \( \pi \) is denoted by \( \text{asc}(\pi) \) (resp., \( \text{des}(\pi) \)). Additionally, \( \pi_i \) is a {\em right-to-left maximum} (resp., {\em right-to-left minimum}) in \( \pi \) if \( \pi_i \) is greater (resp., smaller) than any element to its right. Note that \( \pi_n \) is always a right-to-left maximum and a right-to-left minimum. The number of right-to-left maxima and right-to-left minima in \( \pi \) are denoted by \( \text{rlmax}(\pi) \) and \( \text{rlmin}(\pi) \), respectively. Similarly, we define left-to-right maxima, left-to-right minima, \( \text{lrmax}(\pi) \), and \( \text{lrmin}(\pi) \). For example, if \( \pi = 34152 \), then \( \text{lrmax}(\pi) = 3 \) and \( \text{lrmin}(\pi) = \text{rlmin}(\pi) = \text{rlmax}(\pi) = \text{asc}(\pi) = \text{des}(\pi) = 2 \).

We are interested in the following g.f.:
$$
F_k(x,p,q,u,v,s,t):=\sum_{n\geq 0}x^n\sum_{\sigma \in S_n(2413,3142,P_k) }p^{\asc(\sigma)}q^{\des(\sigma)}u^{\lmax(\sigma)}v^{\rmax(\sigma)}s^{\lmin(\sigma)}t^{\rmin(\sigma)},
$$
for $k\geq 0$. Our studies extend the results of Chen, Kitaev, and Zhang~\cite{DAM2024CZS} on separable permutations. Moreover, they also extend the results of Han and Kitaev~\cite{DMTCS2024HS} on pattern-avoiding permutations. For brevity, throughout the paper, we often omit the variables in the arguments that are set to 1. For example, $F_k(x,p,q,1,v,1,1)$ can be denoted by $F_k(x,p,q,v)$. This will not cause any confusion.

Note that if $k=1$, the only separable permutation avoiding $P_1$  is the empty permutation $\varepsilon$. It follows that $F_1(x,p,q,u,v,s,t)=1$. Moreover, if $k=2$, only the decreasing permutation $n(n-1)\cdots 1$ belongs to $S_n(2413,3142,P_2)$. Hence,
\begin{align}
F_2(x, p, q, u, v, s, t)
= 1+\sum_{n=1}q^{n-1}uv^{n}s^{n}tx^{n}
=1+\frac{uvstx}{1 - qvsx}.
\end{align}
The cases of $k=3$ and $k=4$ will be given by Theorems~\ref{F3-thm} and~\ref{dist-general-xyuv-thm}, respectively, along with a system of functional equations in Section~\ref{most-general-sec} giving $F_k(x,p,q,u,v,s,t)$ for any $k\geq 2$. Notably, our formula for $F_4(x,p,q,u,v,s,t)$ is a rational function, with the sum in the numerator (resp., denominator) containing  100 (resp., 19) monomials.

This paper is organized as follows. In Section~\ref{F3-F4-sec}, we derive formulas for $F_3(x,p,q,u,v,s,t)$ and $F_4(x,p,q,u,v,s,t)$. In Section~\ref{most-general-sec}, we not only provide an expression for $F_k(x,p,q,u,v,s,t)$ for any $k\geq 2$ in terms of a system of functional equations but also demonstrate how to use it to re-derive $F_4(x,p,q,u,v,s,t)$. Additionally, we consider a special general case of $F_k(x,u,v)$, as the joint distribution of left-to-right maxima and right-to-left maxima has attracted significant attention in the  literature for various combinatorial structures; e.g., see \cite{CarlitzScoville,DAM2024CZS, DMTCS2024HS, HanKitZha}. This case involves a much simpler system of functional equations, and we illustrate its application by deriving a formula for $F_3(x,u,v)$.
	
\section{The generating functions $F_3(x,p,q,u,v,s,t)$ and $F_4(x,p,q,u,v,s,t)$}\label{F3-F4-sec}

The following theorem is used in the proof of Theorem~\ref{dist-general-xyuv-thm}.

\begin{thm}\label{F3-thm}
The g.f.\ $F_3(x,p,q,u,v,s,t)$ is given by
\begin{align}\label{F3}
  \frac{1+q^2v^2sx^2+vtusx(1+ptux)-qvx\left(1+pus^2x^2vt(-1+t)
   (-1+u)+s(1+px+vtux)\right)}{1+q^2v^2sx^2-qvx(1+s+psx)}. 
\end{align}
The initial terms of $F_3(x,p,q,u,v,s,t)$ are $1+uvstx+(pu^2vst^2+quv^2s^2t)x^2+(pqu^2v^2st^2+pqu^2v^2s^2t+pquv^2s^2t^2+q^2uv^3s^3t)x^3+\cdots$.
\end{thm}

\proof
By \cite[Thm 2.4]{DMTCS2024HS}, the g.f.\
\begin{align*}
   F_{(123,132)}(x,p,q,u,v,s,t):=
   \sum_{n\geq 0}x^n\sum_{\sigma \in S_n(123,132) }p^{\asc(\sigma)}q^{\des(\sigma)}u^{\lmax(\sigma)}v^{\rmax(\sigma)}s^{\lmin(\sigma)}t^{\rmin(\sigma)}
\end{align*}
is given by
\begin{align*}
   \frac{1+q^2s^2vx^2+stuvx(1+ptux)-qsx(1+puv^2x^2st(-1+t)(-1+u)+v(1+px+stux))}{1+q^2s^2vx^2-qsx(1+v+pvx)}.
\end{align*}
Given any permutation $\pi$ in $S_n(123,132) $, after taking reverse and then complement operations on $\pi$, the new permutation $\pi$ is in $S_n(123,213) $. Moreover, the six statistics (asc, des, lmax, rmax, lmin, rmin) are changed to (asc, des, rmin, lmin, rmax, lmax),
as illustrated in Table~\ref{tab-transf}.

\begin{table}[htbp]
\centering
\renewcommand{\arraystretch}{0}
				\begin{tabular}
{>{\centering\arraybackslash}p{2.5cm}
>{\centering\arraybackslash}p{1cm}
>{\centering\arraybackslash}p{1cm}
>{\centering\arraybackslash}p{1cm}
>{\centering\arraybackslash}p{1cm}>
{\centering\arraybackslash}p{1cm}>
{\centering\arraybackslash}p{1cm}}
		\toprule
		& $p$ & $q$ & $u$ & $v$ & $s$ & $t$ \\
\midrule
123,132 & asc & des & lmax & rmax & lmin & rmin \\
\tikz[baseline]{\draw[-, thick](0,0) -- (0,-0.3);}\\
	{\footnotesize (reverse)} & & & & & &\\
\tikz[baseline]{\draw[->, thick](0,0) -- (0,-0.3);}\\
321,231 & des & asc & rmax & lmax & rmin & lmin \\
\tikz[baseline]{\draw[-, thick](0,0) -- (0,-0.3);}\\
	{\footnotesize (complement)} & & & & & &\\
	\tikz[baseline]{\draw[->, thick](0,0) -- (0,-0.3);}\\
	123,213 & asc & des & rmin & lmin & rmax & lmax \\
	\bottomrule
	\end{tabular}
\caption{The transformation of the six statistics from $S_n(123, 132)$ to $S_n(123, 213)$ under the reverse and  complement operations.}
	\label{tab-transf}
	\end{table}

\noindent Hence, for
\begin{align*}
   F_{(123,213)}(x,p,q,u,v,s,t):=
   \sum_{n\geq 0}x^n\sum_{\sigma \in S_n(123,213) }p^{\asc(\sigma)}q^{\des(\sigma)}u^{\lmax(\sigma)}v^{\rmax(\sigma)}s^{\lmin(\sigma)}t^{\rmin(\sigma)},
\end{align*}
we have
\begin{align}\label{gf-all-5-lem}
F_{(123,213)}(x,p,q,u,v,s,t)
=
F_{(123,132)}(x,p,q,t,s,v,u).
\end{align}
Note that for any permutation $\pi \in S_n(P_3)=S_n(123,213)$, $\pi$ automatically avoids 2413 and 3142, so $S_n(123,213)=S_n(2413, 3142, P_3)$. Therefore, we have
\begin{align*}
  F_3(x,p,q,u,v,s,t)=F_{(123,213)}(x,p,q,u,v,s,t)=F_{(123,132)}(x,p,q,t,s,v,u),
  \end{align*}
which completes the proof.
\qed

\begin{cor}
The number of permutations in $S_n(2413,3142,P_3)$ is given by the g.f.
\begin{align*}
	F_3(x) &=
	\frac{1-x}{1-2x},
\end{align*}
where the initial terms are $1 + x +2x^2 + 4 x^3 + 8x^4 +16 x^5 +32 x^6 + 64 x^7 +128x^8 +256x^9 +512x^{10} + \cdots$.
\end{cor}

\proof
Setting $p=q=u=v=s=t=1$ in  the formula  for $F_3(x,p,q,u,v,s,t)$ in Theorem
	\ref{F3-thm}, we get the desired result.
\qed

We can now prove the following theorem.

\begin{thm}\label{dist-general-xyuv-thm}
For $S_n(2413,3142,P_4)$, we have that the g.f.\
\begin{align}
\label{gf-all-stat-eqn}
F_4(x,p,q,u,v,s,t)= \frac{F_a}{F_b},
\end{align}
where
\begin{align*}
F_b =&
1 + q^4 s v^3 x^4 (v - p s x)  - q v x (3 + s + p x + p s x + p^2 s x^2)\\
&- q^3 v^2 x^3 \left( v - p s^2 x (2 + p x ) + s v (3 + 2 p x) \right)\\
&+ q^2 v x^2 \Big( - p s^2 x + v \left(3 + p x + s (3 + 3 p x + p^2 x^2) \right) \Big),
\end{align*}
\begin{align*}
F_a=&
1 + q^4 s v^3 x^4 \left(v - p s x + p s^2 t u (-1 + t + u - t u) v x^2 \right)\\
&+ s t u v x (1 + p t u x + p^2 t^2 u^2 x^2)\\
&+ q^2 v x^2 \Big[- p s^2 x
+ s t u v^2 x \Big( 3 + p \big( 1 + 2 s (-1 + t) (-1 + u)+ t u \big) x \\
&+ p^2 (-1 + t ) (-1 + u) \big( t u + s (1 + t + u + 2 t u) \big) x^2 \\
&+ p^3 s (-1 + t) t (-1 + u) u x^3 \Big) + v \big(3 + p x - p s^2 t^2 u^2 x^2 (2 + p x) \\
&+ s (3 +3 p x + p^2 x^2) - p s^3 (-1 + t ) t (-1 + u ) u x^2 (1 + p^2 t u x^2 ) \big)\Big]\\
&- q^3 v^2 x^3 \Big(-p s^2 x ( 2 + p x) + s t u v^2 x \big( 1 + p s (-1 + t)(-1 + u) x \\
&+p^2 s (-1 + t) t (-1 + u) u x^2\big) + v \big(1- p s^2 t^2 u^2 x^2 +s (3 + 2 p x )\\
&- p s^3 (-1 + t) t (-1 + u) u x^2 (2 + p x + p^2 t u x^2 ) \big) \Big)\\
&-q v x \Big[3 + p x + s \Big(1 + 3 t u v x + p x ( 1 + t u v x + 2 t^2 u^2 v x )\\
&+ p^2 x^2 \big(1 - t^2 (-1 + u) u^2 v x + t^3 u^2 (-1 + 2u) v x \big) \Big) \\
&+ p s^2 t u x^2 \Big( p t^2 (-1 + u) v x (1 + u + p u x) - (-1 + u) v \big(1 + p (1 + u) x\big)\\
&- t \big( v + p^2 u^2 v x^2- u (-1 + v + p^2 v x^2 ) \big) \Big) \Big].
\end{align*}
\end{thm}

\proof
Our proof relies on Stankova's decomposition of separable permutations~\cite{Stankova}, which is presented in Figure~\ref{fig:stankova}.

When $n=0$, the empty permutation $\varepsilon$ contributes 1 to  $F_4(x,p,q,u,v,s,t)$.
When  $n \geq 1$,   note that in order to avoid $P_4$,  the largest element $n$ of any  permutation $\sigma \in S_n(2413,3142,P_4)$ can only appear in the first three positions of $\sigma=\sigma_1\cdots\sigma_n$.
Hence,   we have the following three disjoint possibilities, where Cases 1--3 cover all permutations in $S_n(2413,3142,P_4)$.

\begin{figure}[htbp]
  \centering
	\begin{tikzpicture}[line width=0.5pt,scale=0.24]
				\coordinate (O) at (0,0);
				\path (10,1)  node {$n$};
				\draw [dashed] (8,0)--++(10,0);
				\fill[black!100] (O)++(10,0) circle(1.5ex);
				\draw (11,-1) rectangle (15,-4);
				\path (13,-2.5)  node {$A$};
				\draw [dashed] (10,0)--++(0,-5);
		\end{tikzpicture}
				\caption{The permutation diagram of $\sigma \in S_n(2413,3142,P_4)$ with $\sigma_1=n$.}
                \label{sepStructure_pattern_4_1}
			\end{figure}
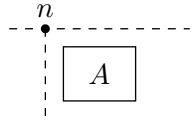

\noindent \textbf{Case 1.} $\sigma_1=n$, as shown in Figure~\ref{sepStructure_pattern_4_1}. Referring to the structure of  a separable permutation in Figure~\ref{sepStructure}, there are two subcases to consider.
\begin{itemize}
\item[(a)] $A=\varepsilon$, then $\sigma=1$ giving the term of $uvstx$.
\item[(b)] $A\neq \varepsilon$, then $\sigma=n  A=1 \ominus A$. In this subcase, we have the term
$$quvsx\left(F_4(x,p,q,1,v,s,t)-1\right)$$
since $\lmax(A)$  does not contribute to $\lmax(\sigma)=1$ and $A$ can be any  permutation in $S_{n-1}(2413,3142,P_4)$, while $n$ contributes an extra descent, one extra left-to-right maximum, one extra right-to-left maximum, and one  extra left-to-right minimum.
\end{itemize}

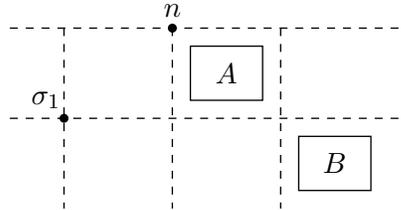
\begin{figure}[htbp]
\centering
 \begin{tikzpicture}[line width=0.5pt,scale=0.24]
		\coordinate (O) at (0,0);
		\path (15,1)  node {$n$};
		\path (8,-4)  node {$\sigma_1$};
		\draw [dashed] (6,0)--++(22,0);
		\draw [dashed] (6,-5)--++(22,0);
		\fill[black!100] (O)++(15,0) circle(1.5ex);
		\fill[black!100] (O)++(9,-5) circle(1.5ex);
		\draw (16,-1) rectangle (20,-4);
		\path (18,-2.5)  node {$A$};
		\draw (22,-6) rectangle (26,-9);
	\path (24,-7.5)  node {$B$};
		\draw [dashed] (15,0)--++(0,-10);
		\draw [dashed] (9,0)--++(0,-10);
		\draw [dashed] (21,0)--++(0,-10);
		\end{tikzpicture}
		\caption{The permutation diagram of $\sigma \in S_n(2413,3142,P_4)$ with $\sigma_2=n$.}\label{sepStructure_pattern_4_2}
			\end{figure}

\noindent \textbf{Case 2.}   $\sigma_2=n$, as shown in Figure~\ref{sepStructure_pattern_4_2}.
Referring to the structure of  a separable permutation in Figure~\ref{sepStructure}, there are four subcases to consider.

\begin{itemize}
\item[(a)] $A=B=\varepsilon$, then  $\sigma=12$ giving the term of  $pu^2vst^2x^2$.
			
\item[(b)]  $A=\varepsilon$ and  $B\neq \varepsilon$,
then $\sigma=(n-1)nB=12 \ominus B$. In this subcase, we have the term
$$pqu^2vsx^2\left(F_4(x,p,q,1,v,s,t)-1\right)$$
since
$\lmax(B)$  does not contribute to $\lmax(\sigma)=1$ and $B$ can be any non-empty permutation in $S_{n-2}(2413,3142,P_4)$,
while the subpermutation $(n-1)n$ contributes an extra ascent, one extra descent, two extra left-to-right maxima, one extra right-to-left maximum, and  one extra left-to-right minimum.

\item[(c)] $A\neq \varepsilon$ and $B=\varepsilon$, then there exists no occurrence of the pattern $P_3$ in $A$, or it would form an occurrence of $P_4$ with the element $\sigma_1$ in $\sigma$.
In this subcase, we have the term
$$pqu^2vstx^2\left(F_3(x,p,q,1,v,1,t)-1\right)$$
since
$\lmax(A)$ and $\lmin(A)$ do not contribute to $\lmin(\sigma)= 1$, $\lmax(\sigma)=2$, and $\lmin(\sigma)=1$, and $A$ can be any non-empty permutation in $S_{n-2}(2413,3142,P_3)$,
while the elements $\sigma_1$ and $n$ contribute an extra ascent, one extra descent, two extra left-to-right maxima, one extra right-to-left maximum,  one extra left-to-right minimum, and  one extra right-to-left minimum.

\item[(d)] $A\neq \varepsilon$ and  $B\neq \varepsilon$,
we have the term
$$
pq^2u^2vsx^2\left(F_3(x,p,q,1,v,1,1)-1\right)
\left(F_4(x,p,q,1,v,s,t)-1\right)$$
since
 \begin{itemize}
 \item[$\bullet$]  the elements $\sigma_1$ and $n$ contribute an extra ascent, one extra descent, two extra left-to-right maxima, one extra right-to-left maximum,  and one extra left-to-right minimum;
 \item[$\bullet$] $A$ does not contribute to the statistics lmax, lmin, and rmin in $\sigma$ and $A$ can be any  permutation in $S_{n-2}(2413,3142,P_3)$. Hence,  we have
the g.f.\ $F_3(x,p,q,1,v,1,1)-1$;
 \item[$\bullet$] $B$ does not contribute to  lmax in $\sigma$ and hence has
g.f.\ $F_4(x,p,q,1,v,s,t)-1$,  and
 \item[$\bullet$] an extra descent is formed between $A$ and $B$.
 \end{itemize}
\end{itemize}

\noindent \textbf{Case 3.} $\sigma_3=n$. We divide Case 3 into  Subcase 3.1 and Subcase 3.2.

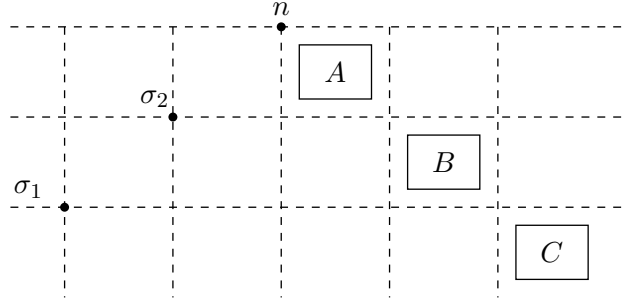
\begin{figure}[htbp]
\centering
\begin{tikzpicture}[line width=0.5pt,scale=0.24]
	\coordinate (O) at (0,0);
						\path (15,1)  node {$n$};
						\path (8,-4)  node {$\sigma_2$};
						\path (1,-9)  node {$\sigma_1$};
						\draw [dashed] (0,0)--++(35,0);
						\draw [dashed] (0,-5)--++(35,0);
						\draw [dashed] (0,-10)--++(35,0);
						\fill[black!100] (O)++(15,0) circle(1.5ex);
						\fill[black!100] (O)++(9,-5) circle(1.5ex);
						\fill[black!100] (O)++(3,-10) circle(1.5ex);
						\draw (16,-1) rectangle (20,-4);
						\path (18,-2.5)  node {$A$};
						\draw (22,-6) rectangle (26,-9);
						\path (24,-7.5)  node {$B$};
						\draw (28,-11) rectangle (32,-14);
						\path (30,-12.5)  node {$C$};
						\draw [dashed] (15,0)--++(0,-15);
						\draw [dashed] (9,0)--++(0,-15);
						\draw [dashed] (3,0)--++(0,-15);
						\draw [dashed] (21,0)--++(0,-15);
						\draw [dashed] (27,0)--++(0,-15);
		\end{tikzpicture}
    \caption{The permutation diagram of $\sigma \in S_n(2413,3142,P_4)$ with $\sigma_3=n$ and $\sigma_1 < \sigma_2$.}
    \label{sepStructure_pattern_4_3_1}	\end{figure}

\noindent
\textbf{Subcase 3.1} 
Suppose $\sigma_1<\sigma_2$ as shown in Figure~\ref{sepStructure_pattern_4_3_1}.
Referring to the structure of  a separable permutation in Figure~\ref{sepStructure}, there are eight subcases to consider.

\begin{itemize}
 \item[(a)] $A=B=C=\varepsilon$,  then  $\sigma=123$ giving the term of  $p^2u^3vst^3x^3$.

\item[(b)]  $A=B=\varepsilon$ and $C\neq \varepsilon$,
then $\sigma= \sigma_1\sigma_2 n C$. In this subcase, we have the term
$$p^2qu^3vsx^3\left(F_4(x,p,q,1,v,s,t)-1\right)$$
since $C$ does not contribute to the statistic lmax  in $\sigma$ and $C$ can be any  permutation in $S_{n-3}(2413,3142,P_4)$,
while the subpermutation $\sigma_1\sigma_2n$ contributes two extra ascents, one extra descent, three extra left-to-right maxima, one extra right-to-left maximum, and  one extra left-to-right minimum.

\item[(c)]
$A=\varepsilon$, $B\neq \varepsilon$, and $C= \varepsilon$, then $\sigma=\sigma_1\sigma_2nB$.
There exists no occurrence of $P_3$ in $B$, since otherwise it would form $P_4$ together with the element $\sigma_1$ in $\sigma$.
In this subcase, we have the term
$$
p^2qu^3vstx^3\left(F_3(x,p,q,1,v,1,t)-1\right)
$$
since
$C$ does not contribute to the statistics lmax  and lmin in $\sigma$, and $C$ can be any  permutation in $S_{n-3}(2413,3142,P_3)$,
while the elements $\sigma_1$, $\sigma_2$ and $n$ contribute two extra ascents, one extra descent, three extra left-to-right maxima, one extra right-to-left maximum,  one extra left-to-right minimum, and  one extra right-to-left minimum.

\item[(d)]
$A=\varepsilon$, $B\neq \varepsilon$, and $C\neq \varepsilon$, then we have the term
$$p^2q^2u^3vsx^3\left(F_3(x,p,q,1,v,1,1)-1\right)\left(F_4(x,p,q,1,v,s,t)-1\right)$$
 since
 \begin{itemize}
 \item[$\bullet$]  the elements $\sigma_1$, $\sigma_2$, and $n$ contribute two extra ascents, one extra descent, three extra left-to-right maxima, one extra right-to-left maximum,  and one extra left-to-right minimum;
 \item[$\bullet$] $B$ does not contribute to the statistics lmax, lmin, and rmin in $\sigma$ and $B$ can be any  permutation in $S_{n-3}(2413,3142,P_3)$. Hence,  we have
the g.f.\ $F_3(x,p,q,1,v,1,1)-1$;
 \item[$\bullet$] $C$ does not contribute to lmax in $\sigma$ and hence has
the g.f.\ $F_4(x,p,q,1,v,s,t)-1$,  and
 \item[$\bullet$] an extra descent is formed between $B$ and $C$.
 \end{itemize}

\item[(e)] $A\neq \varepsilon$ and $B=C=\varepsilon$,
then we have the term
$$p^2qu^3vst^2x^3\cdot\frac{vtx}{1-qvx}
=
\frac{p^2qu^3v^2st^3x^4}{1-qvx}$$
since
\begin{itemize}
\item[$\bullet$] the elements in $A$ must be arranged in descending order, otherwise they would form $P_4$ together with the elements $\sigma_1$ and $\sigma_2$ in  $\sigma$. Hence, $A$  has the g.f.
$$\sum_{i=1}q^{i-1}v^itx^i=\frac{vtx}{1-qvx};$$
\item[$\bullet$] the elements $\sigma_1$, $\sigma_2$ and $n$ contribute two extra ascents, one extra descent, three extra left-to-right maxima, one extra right-to-left maximum,   one extra left-to-right minimum, and  two extra right-to-left minima.
\end{itemize}

\item[(f)] $A\neq \varepsilon$, $B=\varepsilon$, and $ C\neq\varepsilon$, then
we have the term
$$\frac{p^2q^2u^3v^2sx^4}{1-qvx}\left(F_4(x,p,q,1,v,s,t)-1\right)$$
since
\begin{itemize}
\item[$\bullet$] the subpermuation $\sigma_1\sigma_2nA$
does not contribute to the statistic rmin in $\sigma$ and has
the g.f.\ $\frac{p^2qu^3v^2sx^4}{1-qvx}$ from Subcase 3.1(f);
\item[$\bullet$]  $C$ does not contribute to the statistic lmax in $\sigma$ and  can be any non-empty separable  permutation avoiding pattern  $P_4$. Hence, $C$ has
the g.f.\ $F_4(x,p,q,1,v,s,t)-1$,  and
 \item[$\bullet$] an extra descent is formed between $\sigma_1\sigma_2nA$ and $C$.
\end{itemize}

\item[(g)] $A\neq \varepsilon$, $B\neq\varepsilon$, and $ C=\varepsilon$,
then we have the term
$$\frac{p^2q^2u^3v^2stx^4}{1-qvx}\left(F_3(x,p,q,1,v,1,t)-1\right)$$
since
\begin{itemize}
\item[$\bullet$] the elements $\sigma_1$, $\sigma_2$ and $n$ contribute two extra ascents, one extra descent, three extra left-to-right maxima, one extra right-to-left maximum,   one extra left-to-right minimum, and  one extra right-to-left minimum;
\item[$\bullet$] $A$ must be a decreasing permutation and does not contribute to the statistics rmin  in $\sigma$. Hence, $A$  has the g.f.\
$\sum_{i=1}q^{i-1}v^ix^i=\frac{vx}{1-qvx}$ from Subcase 3.1(e);
\item[$\bullet$] $B$ does not contribute to  lmax and lmin in $\sigma$ and  can be any non-empty separable  permutation avoiding $P_3$. Hence, $C$ has the g.f.\ $F_3(x,p,q,1,v,1,t)-1$,  and
\item[$\bullet$] an extra descent is formed between $A$ and $B$.
\end{itemize}

\item[(h)] $A\neq \varepsilon$, $B\neq\varepsilon$, and $ C\neq \varepsilon$, then
we have the term
$$\frac{p^2q^3u^3v^2sx^4}{1-qvx}\left(F_3(x,p,q,1,v,1,1)-1\right)\left(F_4(x,p,q,1,v,s,t)-1\right)$$
since
\begin{itemize}
\item[$\bullet$] the subpermuation $\sigma_1\sigma_2nAB$
does not contribute to the statistic rmin in $\sigma$ and has
the g.f.\ $\frac{p^2q^2u^3v^2sx^4}{1-qvx}\left(F_3(x,p,q,1,v,1,1)-1\right)$ from Subcase 3.1(g);
\item[$\bullet$]  $C$ does not contribute to the statistic lmax in $\sigma$ and  can be any non-empty separable  permutation avoiding $P_4$. Hence, $C$ has
the g.f.\ $F_4(x,p,q,1,v,s,t)-1$,  and
\item[$\bullet$] an extra descent is formed between $\sigma_1\sigma_2nAB$ and $C$.
\end{itemize}

\end{itemize}

\begin{figure}[htbp]
	\centering
	\begin{tikzpicture}[line width=0.5pt,scale=0.24]
						\coordinate (O) at (0,0);
						\path (15,1)  node {$n$};
						\path (8,-9)  node {$\sigma_2$};
						\path (1,-4)  node {$\sigma_1$};
						\draw [dashed] (0,0)--++(35,0);
						\draw [dashed] (0,-5)--++(35,0);
						\draw [dashed] (0,-10)--++(35,0);
						\fill[black!100] (O)++(15,0) circle(1.5ex);
						\fill[black!100] (O)++(9,-10) circle(1.5ex);
						\fill[black!100] (O)++(3,-5) circle(1.5ex);
						\draw (16,-1) rectangle (20,-4);
						\path (18,-2.5)  node {$A$};
						\draw (28,-11) rectangle (32,-14);
						\path (30,-12.5)  node {$B$};
						\draw [dashed] (15,0)--++(0,-15);
						\draw [dashed] (9,0)--++(0,-15);
						\draw [dashed] (3,0)--++(0,-15);
						\draw [dashed] (21,0)--++(0,-15);
						\draw [dashed] (27,0)--++(0,-15);
					\end{tikzpicture}
			\caption{The permutation diagram of $\sigma \in S_n(2413,3142,P_4)$ with $\sigma_3=n$ and $\sigma_1 > \sigma_2$.}
        \label{sepStructure_pattern_4_3_2}
			\end{figure}
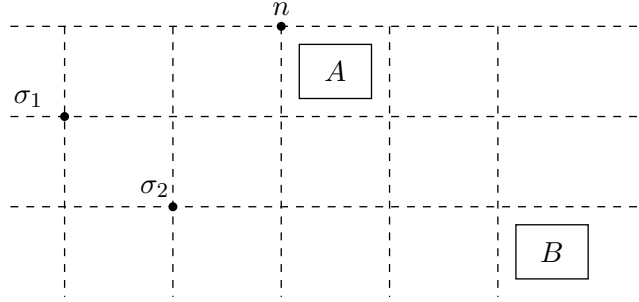

\noindent
\textbf{Subcase 3.2}
Suppose $\sigma_1>\sigma_2$ as shown in Figure~\ref{sepStructure_pattern_4_3_2}.
Note that there exists no  element between $A$ and $B$ which is larger than $\sigma_2$ and smaller than $\sigma_1$, otherwise this element, together with elements $\sigma_1$, $\sigma_2$, and $n$, would form a $3142$ pattern,  a contradiction.
 Referring to the structure of  a separable permutation in Figure~\ref{sepStructure}, there are four subcases to consider.

\begin{itemize}

\item[(a)] When $A=B=\varepsilon$,
 then  $\sigma=213$ giving the term of  $pqu^2vs^2t^2x^3$.

\item[(b)] When $A=\varepsilon$ and $B \neq \varepsilon$,
we have the term
$$pq^2u^2vs^2x^3\left(F_4(x,p,q,1,v,s,t)-1\right)$$
since  $B$ does not contribute to the statistic lmax   in $\sigma$ and  can be any  permutation in $S_{n-3}(2413,3142,P_4)$,
while the elements $\sigma_1$, $\sigma_2$ and $n$ contribute an extra ascent, two extra descents, two extra left-to-right maxima, one extra right-to-left maximum,   and two extra left-to-right minima.

\item[(c)]
When $A \neq \varepsilon$ and $B = \varepsilon$,
we have the term
$$\frac{pq^2u^2v^2s^2t^2x^4}{1-qvx}$$
since
the elements in $A$ is arranged in descending order and $A$ does not contribute to the statistic lmax and lmin in $\sigma$. Hence, $A$  has the g.f.\
$\sum_{i=1}q^{i-1}v^itx^i=\frac{vtx}{1-qvx}$
from Subcase 3.1(e),
while the elements $\sigma_1$, $\sigma_2$ and $n$ contribute an extra ascent, two extra descents, two extra left-to-right maxima, one extra right-to-left maxima,   two extra left-to-right minima, and  one extra right-to-left minimum.

\item[(d)]
When $A \neq \varepsilon$ and $B \neq \varepsilon$,
we have the term
 $$\frac{pq^3u^2v^2s^2x^4}{1-qvx}\left(F_4(x,p,q,1,v,s,t)-1\right)$$
 since
 \begin{itemize}
 \item[$\bullet$] the elements $\sigma_1$, $\sigma_2$ and $n$ contribute one extra ascent, two extra descent, two extra left-to-right maxima, one extra right-to-left maximum,   and two extra left-to-right minima,
\item[$\bullet$] $A$ must is arranged in descending order and does not contribute to the statistics lmax, lmin, and rmin  in $\sigma$. Hence, $A$  has the g.f.\
$\sum_{i=1}q^{i-1}v^ix^i=\frac{vx}{1-qvx}$ from Subcase 3.1(e),
\item[$\bullet$] $B$ does not contribute to the statistics lmax in $\sigma$ and  can be any non-empty separable  permutation avoiding pattern  $P_4$. Hence, $B$ has the g.f.\ $F_4(x,p,q,1,v,s,t)-1$,  and
\item[$\bullet$] an extra descent is formed between $A$ and $B$.
\end{itemize}
		
\end{itemize}

Summarising all the subcases in Cases 1--3, we obtain
\begin{align}\label{F_4_gf_or}
F_4(x,p,q,u,v,s,t)=&
1+a_4(x,p,q,u,v,s,t)+b_4(x,p,q,u,v,s,t)
\nonumber \\
&~~~~~~~~+c_4(x,p,q,u,v,s,t) +d_4(x,p,q,u,v,s,t),
\end{align}

where
\begin{align*}
a_4(x,p,q,u,v,s,t):=&uvstx+quvsx\left(F_4(x,p,q,1,v,s,t)-1\right),
\\
b_4(x,p,q,u,v,s,t):=&pu^2vst^2x^2+pqu^2vsx^2\left(F_4(x,p,q,1,v,s,t)-1\right)\\
&+pqu^2vstx^2\left(F_3(x,p,q,1,v,1,t)-1\right)\\
&+pq^2u^2vsx^2\left(F_3(x,p,q,1,v,1,1)-1\right)
\left(F_4(x,p,q,1,v,s,t)-1\right), \\
c_4(x,p,q,u,v,s,t):=&p^2u^3vst^3x^3+p^2qu^3vsx^3\left(F_4(x,p,q,1,v,s,t)-1\right)\\
&+p^2qu^3vstx^3\left(F_3(x,p,q,1,v,1,t)-1\right)
\\
&+p^2q^2u^3vsx^3\left(F_3(x,p,q,1,v,1,1)-1\right)\left(F_4(x,p,q,1,v,s,t)-1\right)\\
&+\frac{p^2qu^3v^2st^3x^4}{1-qvx}+\frac{p^2q^2u^3v^2sx^4}{1-qvx}\left(F_4(x,p,q,1,v,s,t)-1\right)\\
&+\frac{p^2q^2u^3v^2stx^4}{1-qvx}\left(F_3(x,p,q,1,v,1,t)-1\right)\\
&+\frac{p^2q^3u^3v^2sx^4}{1-qvx}\left(F_3(x,p,q,1,v,1,1)-1\right)\left(F_4(x,p,q,1,v,s,t)-1\right),
\end{align*}

\begin{align*}
d_4(x,p,q,u,v,s,t):=&
pqu^2vs^2t^2x^3+pq^2u^2vs^2x^3\left(F_4(x,p,q,1,v,s,t)-1\right)\\
&+\frac{pq^2u^2v^2s^2t^2x^4}{1-qvx}
+
\frac{pq^3u^2v^2s^2x^4}{1-qvx}\left(F_4(x,p,q,1,v,s,t)-1\right).\end{align*}

\noindent Letting $u=1$ in \eqref{F_4_gf_or},
we get
\begin{align}\label{F_4_gf_or_u}
F_4(x,p,q,v,s,t)=&
1+a_4(x,p,q,v,s,t)+b_4(x,p,q,v,s,t)
\nonumber \\
&~~~~~~~~+c_4(x,p,q,v,s,t) +d_4(x,p,q,v,s,t),
\end{align}
where\\[-5mm]
\begin{align*}
a_4(x,p,q,v,s,t):=&vstx+qvsx\left(F_4(x,p,q,v,s,t)-1\right),
\\
b_4(x,p,q,v,s,t):=&pvst^2x^2+pqvsx^2\left(F_4(x,p,q,v,s,t)-1\right)\\
&+pqvstx^2\left(F_3(x,p,q,v,1,t)-1\right)\\
&+pq^2vsx^2\left(F_3(x,p,q,v,1,1)-1\right)
\left(F_4(x,p,q,v,s,t)-1\right), \\
c_4(x,p,q,v,s,t):=&p^2vst^3x^3+p^2qvsx^3\left(F_4(x,p,q,v,s,t)-1\right)\\
&+p^2qvstx^3\left(F_3(x,p,q,v,1,t)-1\right)
\\
&+p^2q^2vsx^3\left(F_3(x,p,q,v,1,1)-1\right)\left(F_4(x,p,q,v,s,t)-1\right)\\
&+\frac{p^2qv^2st^3x^4}{1-qvx}+\frac{p^2q^2v^2sx^4}{1-qvx}\left(F_4(x,p,q,v,s,t)-1\right)\\
&+\frac{p^2q^2v^2stx^4}{1-qvx}\left(F_3(x,p,q,v,1,t)-1\right)\\
&+\frac{p^2q^3v^2sx^4}{1-qvx}\left(F_3(x,p,q,v,1,1)-1\right)\left(F_4(x,p,q,v,s,t)-1\right), 
\end{align*}
\begin{align*}
d_4(x,p,q,v,s,t):=&
pqvs^2t^2x^3+pq^2vs^2x^3\left(F_4(x,p,q,v,s,t)-1\right)\\
&+\frac{pq^2v^2s^2t^2x^4}{1-qvx}
+
\frac{pq^3v^2s^2x^4}{1-qvx}\left(F_4(x,p,q,v,s,t)-1\right).  
\end{align*}
Combining with Theorem~\ref{F3-thm}, we obtain
\begin{align}\label{gf-all-stat-eqn-F4-u}
F_4(x,p,q,v,s,t)= \frac{F_c}{F_d},
\end{align}
where\\[-7mm]
\begin{align*}
F_c=&
1 + q^4 s v^3 x^4 (v - p s x) + s t v x (1 + p t x + p^2 t^2 x^2) \\
&+ q^2 v x^2 \Big(- p s^2 x + s t v^2 x (3 + px +pt x) +  v \left(3 + p x - p s^2 t^2  x^2 (2 + p x) + s (3 + 3 p x +p^2 x^2)\right)\Big) \\
&-
q^3 v^2 x^3 \Big(v + s t v^2 x - p s^2 t^2 v x^2 - p s^2 x (2 + p x) + s v (3 + 2 p x)\Big)
 \\
&-
q v x \Big(3 + p x - p s^2 t^2 x^2 + s \left(1 + 3 t v x + p x (1 + t v x + 2 t^2 v x) + p^2 x^2 (1 + t^3 v x )\right)\Big),
\end{align*}
\vspace{-0.5cm}
\begin{align*}
F_d=&
1 + q^4 s v^3 x^4 (v - p s x) - q v x (3 + s + p x + p s x + p^2 s x^2) \\
&+
q^2 v x^2 \Big( - p s^2 x +  v  \left(3 + p x + s (3 + 3 p x + p^2 x^2 )\right )\Big)\\
&-
q^3  v^2 x^3 \Big(v - p s^2 x (2 + p x) + s v (3 +2 p x) \Big).
\end{align*}
Substituting \eqref{gf-all-stat-eqn-F4-u} in \eqref{F_4_gf_or}, we get the expression of
$F_4(x,p,q,u,v,s,t).$
\qed


\begin{cor}
The number of permutations in $S_n(2413,3142,P_4)$ is given by the g.f.
\begin{align*}
F_4(x) &= \frac{1-3x+2x^2}{1-4x+4x^2-2x^3+2x^4},
\end{align*}
where the initial terms are $1 + x +2x^2 + 6 x^3 + 16x^3 +42 x^5 +112 x^6 + 300 x^7 +804x^8 +2156x^9 +5784x^{10} + \cdots$.
\end{cor}

\proof
Setting $p=q=u=v=s=t=1$ in  the formula  for $F_4(x,p,q,u,v,s,t)$ in Theorem \ref{dist-general-xyuv-thm}, we get the desired result.
\qed

\begin{rem}
	\textnormal{
We note that both $F_3(x,u,t)$ and $F_4(x,u,t)$ are symmetric in $u$ and $t$.
Let $p = q = v = s = 1$, and set one of the variables $u$ or $t$ equal to $1$ in \eqref{F3} and \eqref{gf-all-stat-eqn}.
This yields the single distributions of $\operatorname{lrmax}$ and $\operatorname{rlmin}$ over $S_n(123,213)$ and  $S_n(2413,3142,P_4)$:}
	\begin{align*}
F_3(x,u) &= \frac{1 - 2x + ux - ux^2 + u^2x^2}{1 - 2x}, \\
F_3(x,t) &= \frac{1 - 2x + tx - tx^2 + t^2x^2}{1 - 2x},\\
F_4(x,u) &= \frac{1 + (-4 + u)x + (4 - 3u + u^2)x^2 + (-2 + 2u - u^2 + u^3)x^3 - (-2 + u^2 + u^3)x^4}{1 - 4x + 4x^2 - 2x^3 + 2x^4}, \\
F_4(x,t) &= \frac{1 + (-4 + t)x + (4 - 3t + t^2)x^2 + (-2 + 2t - t^2 + t^3)x^3 - (-2 + t^2 + t^3)x^4}{1 - 4x + 4x^2 - 2x^3 + 2x^4}.
\end{align*}
\textnormal{These distributions coincide because the patterns 123 and 213, as well as 2413, 3142, and $P_4$, are invariant under the group-theoretic inverse operation, which exchanges left-to-right maxima with right-to-left minima. For the same reason, $F_k(x,u,t)$ is symmetric in $u$ and $t$ for any $k\geq 2$.}
\end{rem}

\section{The generating function $F_k(x,p,q,u,v,s,t)$ for $k\geq 2$}\label{most-general-sec}

In this section, we derive a system of functional equations for $F_k(x,p, q, u,v,s,t)$
for any $k \geq 2$. To do this, for any $k \geq 2 $  and $0\leq i \leq k-2$, we define
\begin{align*}
f_{k,i}(x,p, q, u,v,s,t) := &\sum_{n \geq 1} x^n \sum_{\pi} p^{\asc(\pi)}q^{\des(\pi)}u^{\lmax(\pi)}v^{\rmax(\pi)}s^{\lmin(\pi)}t^{\rmin(\pi)},\\
g_{k,i}(x,p, q, u,v,s,t) :=& \sum_{n \geq 2} x^n \sum_{\sigma} p^{\asc(\sigma)}q^{\des(\sigma)}u^{\lmax(\sigma)}v^{\rmax(\sigma)}s^{\lmin(\sigma)}t^{\rmin(\sigma)}
\end{align*}
where the sums are over all $\pi$ (resp., $\sigma$)  in $S_n(2413, 3412, P_k)$ with $i$ elements to the left of $n$
and  with the element $n-1$  (if $n-1$ exists) to the right (resp., left) of $n$.
For convenience, we will use $f_{k,i}$ (resp., $g_{k,i}$) to represent $f_{k,i}(x, p, q, u, v, s, t)$ (resp., $g_{k,i}(x, p, q, u, v, s, t)$).
Moreover, let $X_n:=X_n(p, q, u, v, s, t)$ represent the g.f.\ for separable permutations of length  $n$ with respect to six statistics: asc, des, lmax, rmax, lmin, and rmin. Equivalently,  $X_n(p, q, u, v, s, t)$ is the coefficient of $x^n$ in $F_0(x,p, q, u, v, s, t)$.

The main result of this section is the following theorem.

\begin{thm}\label{thm 4}
For any $k\geq 2$,
the g.f.\ $F_k:=F_k(x, p,q,u, v,s,t)$ is given by the system
\begin{align}\label{in_thm_Z-eq-1}
F_k = &1 + \sum\limits_{i=0}^{k-2} \Big(f_{k,i}(x, p,q,u, v,s,t)+ g_{k,i}(x, p,q,u, v,s,t)\Big),
\\
f_{k,i}=& pquvx^{i+1} X_i(p,q,u,1,s,t) \left(F_{k-i}(x,p,q,1,v,1,t) - 1\right) \nonumber  \\
&+
pq^2uvx^{i+1} X_i(p,q,u,1,s,1) \left(F_{k-i}(x,p,q,1,v,1,1) - 1\right)
\left(F_{k}(x,p,q,1,v,s,t)-1\right)
\nonumber\\
&+
pq
\left(F_{k-i}(x,p,q,1,v,1,1) - 1\right)
\sum_{j=1}^{i-1}
X_j(p,q,u,1,1,1) x^j   f_{k,i-j},
\label{eq:fki_term31}
\end{align}

\vspace{-0.7cm}

\begin{align} \label{in_thm_Z-eq-3}
g_{k,i} =& puvtx^{i+1} X_i(p,q,u,1,s,t) + pquvx^{i+1} X_i(p,q,u,1,s,1)  \left(F_k(x,p,q,1,v,s,t) - 1\right) \nonumber \\
&+
p^2quvx^{i+1} \sum_{j=1}^{i-1}
  X_j(p,q,u,1,1,1) X_{i-j}(p,q,u,1,s,t)
  \left(F_{k-i+j}(x,p,q,1,v,1,t) - 1\right) \nonumber\\
  &+
  pq\sum_{j=1}^{i-1} X_{j}(p,q,u,1,1,1) x^j \cdot
\left(
F_{k-i+j}(x,p,q,1,v,1,1) - 1
\right) \nonumber \\
&~~~~~~~\cdot
\left(
g_{k,i-j}
-puvtx^{i-j+1}  X_{i-j}(p,q,u,1,s,t)
\right)
\end{align}
for $1\leq i \leq k-2$,
where the initial conditions are
\begin{align}\label{in_thm_Z-eq-4}
\begin{cases}
f_{k,0}(x, p, q, u, v, s, t)& = uvstx+quvsx\left(F_k(x, p, q, 1, v, s, t)-1\right), \vspace{2mm} \\
g_{k,0}(x, p,q, u, v,s,t)&= 0, \vspace{2mm}\\
f_{2,0}(x, p, q, u, v, s, t) &= \frac{uvstx}{1 - qvsx}.
\end{cases}
\end{align}
\end{thm}

\proof
Similar to the proof of Theorem~\ref{dist-general-xyuv-thm}, we make use of Stankova's decomposition of separable permutations~\cite{Stankova}, as illustrated in Figure~\ref{fig:stankova}.

When $k=2$, the only permutation in $S_n(2413, 3412, P_2)$ is $\pi=n(n-1)\cdots 1$.
Clearly, there exists no element to the left of $n$ in $\pi$, hence the only value for $i$ is 0.
It follows that
\begin{align}\label{thmz-eq-re-0}
g_{2,0}(x, p, q, u, v, s, t) = 0,
~
\mathrm{and}
~
f_{2,0}(x, p, q, u, v, s, t)= \sum_{n\geq 1}q^{n-1}uv^ns^ntx^n= \frac{uvstx}{1 - qvsx}.
\end{align}
Moreover, we have
\begin{align*}
F_2(x, p, q, u, v, s, t) &= 1+f_{2,0}(x, p, q, u, v, s, t) +g_{2,0}(x, p, q, u, v, s, t)
=1+\frac{uvstx}{1 - qvsx},
\end{align*}
which is consistent with  \eqref{in_thm_Z-eq-1}--\eqref{in_thm_Z-eq-3} and the initial conditions \eqref{in_thm_Z-eq-4}.

When $k\geq 3$,
it is clear that
\begin{align}\label{thmz-eq-re-1}
F_k(x, p, q, u, v, s, t)
&=1+\sum_{i=0}^{k-2} f_{k,i}(x, p, q, u, v, s, t)+\sum_{i=0}^{k-2} g_{k,i}(x, p, q, u, v, s, t).
\end{align}

\noindent  If $1 \leq i \leq k-2$,  let $\pi$ be any permutation  in $S_n(2413, 3412, P_k)$ with $i$ elements to the left of $n$. Since $i\geq 1$, it is obvious that $n\geq 2$ and $n-1$ does exist in $\pi$. We have two cases depending on whether the element $n-1$ is to the right or to the left of $n$ in $\pi$.

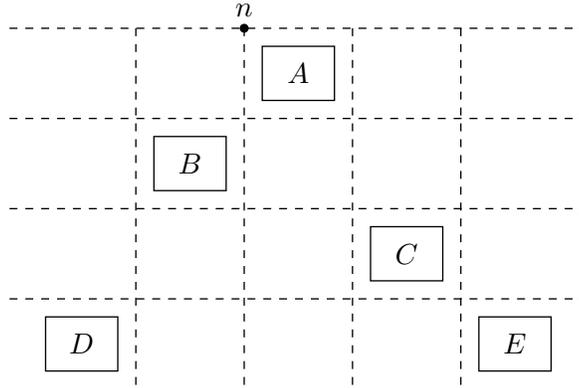
\begin{figure}[htbp]
	\centering
	\begin{tikzpicture}[line width=0.5pt,scale=0.24]
			\coordinate (O) at (0,0);
			\path (15,1)  node {$n$};
			\draw [dashed] (2,0)--++(32,0);
	\draw [dashed] (2,-5)--++(32,0);
						\draw [dashed] (2,-10)--++(32,0);
                            \draw [dashed] (2,-15)--++(32,0);
						\fill[black!100] (O)++(15,0) circle(1.5ex);
						\draw (16,-1) rectangle (20,-4);
						\path (18,-2.5)  node {$A$};
						\draw (10,-6) rectangle (14,-9);
						\path (12,-7.5)  node {$B$};
						\draw (22,-11) rectangle (26,-14);
						\path (24,-12.5)  node {$C$};
                            \draw (28,-16) rectangle (32,-19);
						\path (30,-17.5)  node {$E$};
                            \draw (4,-16) rectangle (8,-19);
						\path (6,-17.5)  node {$D$};
						\draw [dashed] (15,0)--++(0,-20);
						\draw [dashed] (9,0)--++(0,-20);
						\draw [dashed] (21,0)--++(0,-20);
						\draw [dashed] (27,0)--++(0,-20);
					\end{tikzpicture}
				\caption{A schematic representation of the permutation diagrams for permutations in Case 1.}
        \label{sepStructure_thm3_f_1}
				\end{figure}

\noindent \textbf{Case 1.} $n-1$ is to the right  of $n$ in $\pi$, as shown in Figure~\ref{sepStructure_thm3_f_1}.
 Note that $A$ is non-empty since $n-1$ is in $A$. From $i\geq 1$, we get that $B$ is also non-empty.
There are two subcases to consider.

\begin{itemize}
\item[(a)]
$B$ contains $i$ elements,  which implies that $D$ is empty and hence $E$ is also empty. If $C$ is empty, then $A$ does not contain the pattern $P_{k-i}$, or combining with the $i$ elements in $B$ will form a pattern $P_k$.  That is a contradiction.
Since  $A$ does not contribute to the statistics lmax and lmin in $\pi$, and  it can be any  non-empty  separable permutation avoiding  $P_{k-i}$, we have the g.f.\
$F_{k-i}(x,p,q,1,v,1,t) - 1$.
Moreover, $B$ does not contribute to the statistic rmax in $\pi$ and it can be
any non-empty separable permutation of length~$i$. The  g.f.\ is then $X_i(p,q,u,1,s,t)x^i$.
Note that the largest element $n$ contributes  an extra ascent, one extra descent, one extra left-to-right maximum, and one extra right-to-left maximum.
Hence, in this subcase, we have the term
\begin{align}\label{thmZ-1}
&pquvx \cdot X_i(p,q,u,1,s,t)x^i  \cdot \left(F_{k-i}(x,p,q,1,v,1,t) - 1\right) \nonumber\\
&~~~~~~~~~~~~~~~=pquvx^{i+1} X_i(p,q,u,1,s,t) \left(F_{k-i}(x,p,q,1,v,1,t) - 1\right).
\end{align}

If $C$ is non-empty, then $C$
is a non-empty separable
permutation avoiding $P_k$  and
does not contribute to the statistic lmax in $\pi$.
Thus, the g.f.\ is
$F_{k}(x,p,q,1,v,s,t)-1$.
Note that the subpermutation
$BnA$  does not contribute to the statistic
rmin  in $\pi$ and an extra descent is formed between $BnA$ and $CE$. Hence, from \eqref{thmZ-1}, in this subcase we have the term
\begin{align}\label{thmZ-2}
q\left( pquvx^{i+1} X_i(p,q,u,1,s,1) \left(F_{k-i}(x,p,q,1,v,1,1) - 1\right)\right)
\left(F_{k}(x,p,q,1,v,s,t)-1\right).
\end{align}

\item[(b)]
$B$ contains $j$ elements,  where $1\leq j \leq i-1$ (if $i \geq 2$). This implies that $C$ and $D$ are both non-empty. Then, the g.f.\  is
\begin{align}\label{thmZ-3}
pq
\left(F_{k-i}(x,p,q,1,v,1,1) - 1\right)
\sum_{j=1}^{i-1}
X_j(p,q,u,1,1,1) x^j  \cdot f_{k,i-j}(x,p,q,u,v,s,t).
\end{align}
Indeed,
\begin{itemize}
\item[$\bullet$] after removing  $A$ and $B$, $n$ together with all the other elements in $\pi$  is still a  separable
permutation avoiding $P_k$ and hence has the g.f.\
$f_{k,i-j}(x,p,q,u,v,s,t)$;
\item[$\bullet$] $A$ does not contribute to the statistics lmax, lmin, and rmin in $\pi$, and  can be any  non-empty  separable permutation avoiding  $P_{k-i}$. Thus,  the g.f.\ is
$F_{k-i}(x,p,q,1,v,1,1) - 1$;
\item[$\bullet$] $B$ does not contribute to the statistics rmax, lmin, and rmin in $\pi$ and can be
any non-empty separable permutation of length $j$. Thus,  the g.f.\ is
$ X_j(p,q,u,1,1,1) x^{j}$;
\item[$\bullet$]  an extra ascent is formed between $D$ and $B$, and
\item[$\bullet$]  an extra descent is formed between $A$ and $C$.
\end{itemize}
\end{itemize}

\noindent
Combining with \eqref{thmZ-1}, \eqref{thmZ-2}, and \eqref{thmZ-3}, we obtain the following expression for $f_{k,i}(x,p,q,u,v,s,t)$:
\begin{align}
f_{k,i}=& pquvx^{i+1} X_i(p,q,u,1,s,t) \left(F_{k-i}(x,p,q,1,v,1,t) - 1\right) \nonumber  \\
&+
pq^2uvx^{i+1} X_i(p,q,u,1,s,1) \left(F_{k-i}(x,p,q,1,v,1,1) - 1\right)
\left(F_{k}(x,p,q,1,v,s,t)-1\right)
\nonumber\\
&+
pq
\left(F_{k-i}(x,p,q,1,v,1,1) - 1\right)
\sum_{j=1}^{i-1}
X_j(p,q,u,1,1,1) x^j   f_{k,i-j}.\nonumber
\end{align}

\begin{figure}[htbp]
     \centering
\begin{tikzpicture}[line width=0.5pt,scale=0.24]
	\coordinate (O) at (0,0);
		\path (15,1)  node {$n$};
		\draw [dashed] (-3,0)--++(36,0);
        \draw [dashed] (-3,-5)--++(36,0);
						\draw [dashed] (-3,-10)--++(36,0);
                            \draw [dashed] (-3,-15)--++(36,0);
                            \draw [dashed] (-3,-20)--++(36,0);
						\fill[black!100] (O)++(15,0) circle(1.5ex);
						\draw (10,-1) rectangle (14,-4);
						\path (12,-2.5)  node {$A$};
						\draw (16,-6) rectangle (20,-9);
						\path (18,-7.5)  node {$B$};
						\draw (4,-11) rectangle (8,-14);
						\path (6,-12.5)  node {$C$};
                            \draw (22,-16) rectangle (26,-19);
						\path (24,-17.5)  node {$D$};
                            \draw (-2,-21) rectangle (2,-24);
						\path (0,-22.5)  node {$E$};
                            \draw (28,-21) rectangle (32,-24);
						\path (30,-22.5)  node {$F$};
						\draw [dashed] (15,0)--++(0,-25);
						\draw [dashed] (9,0)--++(0,-25);
						\draw [dashed] (3,0)--++(0,-25);
						\draw [dashed] (21,0)--++(0,-25);
						\draw [dashed] (27,0)--++(0,-25);

					\end{tikzpicture}
					\caption{A schematic representation of the permutation diagrams for permutations in Case 2.}
        \label{sepStructure_thm3_g_1}
				\end{figure}
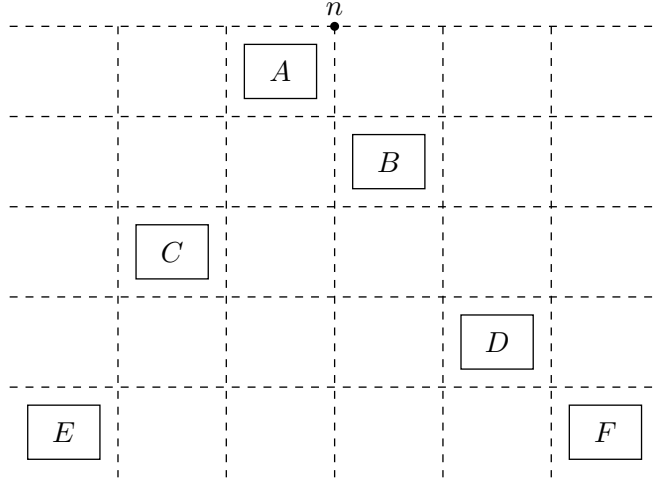

\noindent \textbf{Case 2.} $n-1$ is to the left  of $n$ in $\pi$, as shown in Figure~\ref{sepStructure_thm3_g_1}.
Note that $A$ is non-empty since $n-1$ belongs to $A$.
There are two subcases to consider.

\begin{itemize}
\item[(a)]
$A$ contains $i$ elements,  which implies that $E\cup C$ and hence $D \cup F$ are both empty. If  $B$ is empty, the g.f. is
\vspace{-0.4cm}
\begin{align}\label{thmz_gi_1}
 puvtx \cdot X_i(p,q,u,1,s,t)x^i
\end{align}
since
$A$ does not contribute to  the statistic rmax in $\pi$ and can be any separable permutation of length $i$,
while  the largest element $n$  contributes an extra ascent, one extra left-to-right maximum, one extra right-to-left maximum, and one extra right-to-left minimum.

If  $B$ is non-empty, then $B$ can be any non-empty separable permutation avoiding
$P_k$. Also, $B$ does not contribute to lmax in $\pi$ and hence has the g.f.\ $F_k(x,p,q,1,v,s,t) - 1$.
Note that the subpermutation  $An$ does not contribute to the statistic rmin in $\pi$
and an extra descent is formed between $An$ and $B$.
Hence, from \eqref{thmz_gi_1}, in this subcase we have the term
\begin{align}\label{thmz_gi_2}
 q \cdot puvx^{i+1} X_i(p,q,u,1,s,1)  \cdot \left(F_k(x,p,q,1,v,s,t) - 1\right).
\end{align}

\item[(b)]
$A$ contains $j$ elements,  where $1\leq j \leq i-1$ (if $i \geq 2$). This implies that $B$ and $C$ are both  non-empty.
If $D$ is empty,  then $E$ and $F$ are both empty.
Note that $C$  can be any non-empty separable permutation of length $i-j$ and  does
not contribute to the statistic  rmax in $\pi$, hence the g.f.\
is $X_{i-j}(p,q,u,1,s,t)x^{i-j}$.
Moreover,
$AnB$ does not contribute to the statistic lmin in $\pi$, and  $B$ can be any non-empty separable permutation avoiding $P_{k-(i-j)}$ with $j$ elements to the left of the largest element in
$AnB$, while
an extra ascent is formed between  $C$
and $A$. Hence, from \eqref{thmz_gi_2}, in this subcase we have the term
\begin{small}
\begin{align}\label{thmz_gi_3}
p \sum_{j=1}^{i-1}
\big(
  pquvx^{j+1} X_j(p,q,u,1,1,1)  \left(F_{k-i+j}(x,p,q,1,v,1,t) - 1\right)
\big)
X_{i-j}(p,q,u,1,s,t)x^{i-j}.
\end{align}
\end{small}

If $D$ is non-empty,
this implies that $B$ and $C$ are both non-empty. Then, the g.f.\ is
\begin{align}\label{thmz_gi_4}
pq\sum_{j=1}^{i-1}
&X_{j}(p,q,u,1,1,1) x^j \cdot
\left(
F_{k-i+j}(x,p,q,1,v,1,1) - 1
\right) \nonumber \\
&~~~~\cdot
\left(
g_{k,i-j} (x,p,q,u,v,s,t)
-
puvtx \cdot X_{i-j}(p,q,u,1,s,t)x^{i-j}
\right)
\end{align}
\vspace{-0.3cm}
since\\[-5mm]
\begin{itemize}
\item[$\bullet$] after removing  $A$ and $B$, $n$ together with all the other elements in $\pi$  is still a  separable
permutations avoiding $P_k$ with $i-j$ elements to the left of the largest element and the second largest element is to the left of the largest element.
Note that $D$ is non-empty.
Hence, combining with \eqref{thmz_gi_1}, the g.f.\ is
\begin{align*}
g_{k,i-j} (x,p,q,u,v,s,t)
-
puvtx \cdot X_{i-j}(p,q,u,1,s,t)x^{i-j}
\end{align*}
where the subtraction corresponds to the case of $D$ being empty.
\item[$\bullet$] $A$ does not contribute to the statistics rmax, lmin, and rmin in $\pi$, and  can be any  non-empty  separable permutation  of length $j$. Thus,  the g.f.\ is
$X_{j}(x,p,q,u,1,1,1) x^j$;
\item[$\bullet$] $B$ does not contribute to lmax, lmin, and rmin in $\pi$, and can be any non-empty separable permutation avoiding $P_{k-(i-j)}$. Thus,  the g.f.\ is
$F_{k-i+j}(x,p,q,1,v,1,1) - 1$;
\item[$\bullet$] an extra ascent is formed between $C$ and $A$, and
\item[$\bullet$] an extra descent is formed between $B$ and $D$.
\end{itemize}
\end{itemize}

\noindent
Combining with \eqref{thmz_gi_1}--\eqref{thmz_gi_4}, we obtain the following expression for $g_{k,i}(x,p,q,u,v,s,t)$:
\begin{align}
g_{k,i} =& puvtx^{i+1} X_i(p,q,u,1,s,t) \nonumber \\
&+
pquvx^{i+1} X_i(p,q,u,1,s,1)  \left(F_k(x,p,q,1,v,s,t) - 1\right) \nonumber \\
&+
p^2quvx^{i+1} \sum_{j=1}^{i-1}
  X_j(p,q,u,1,1,1) X_{i-j}(p,q,u,1,s,t)
  \left(F_{k-i+j}(x,p,q,1,v,1,t) - 1\right) \nonumber\\
  &+
  pq\sum_{j=1}^{i-1} X_{j}(p,q,u,1,1,1) x^j \cdot
\left(
F_{k-i+j}(x,p,q,1,v,1,1) - 1
\right) \nonumber \\
&~~~~~~~\cdot
\left(
g_{k,i-j} (x,p,q,u,v,s,t)
-
puvtx^{i-j+1}  X_{i-j}(p,q,u,1,s,t)
\right).\nonumber
\end{align}

It is clear that $X_j(u)$ is the coefficient of $x^j$ in formula (\ref{F(x,u)}) giving the distribution of $\lmax$ on separable permutations. Hence, we can regard $X_j(u)$ as a known function.

\noindent If $i = 0$, let $\pi=\pi_1\cdots\pi_n$ be any non-empty separable permutation avoiding $P_k$ with  $\pi_1=n$. Note that $n-1$ (if  $n-1$ exists ) is to the right of $n$. We have $g_{k,0}(x, p, q, u, v, s, t) = 0$ and
\begin{align}\label{thm3-eq-re-4}
f_{k,0}(x, p, q, u, v, s, t) = uvstx+quvsx\left(F_k(x, p, q, 1, v, s, t)-1\right)
\end{align}
since the subpermutation $\pi_2 \pi_3\cdots\pi_n$ (if $n\geq 2$) in $\pi$
does not contribute to the statistic lmax and can be any permutation in $S_{n-1}(2413, 3142, P_k)$, while the element $n$ contributes an descent, one extra left-to-right maximum, one extra
right-to-left maximum, and one extra left-to-right minimum in $\pi$.

Fixing an integer $k\geq 3$, we suppose that the expressions of $F_{j}(x, p, q, u, v, s, t)$ are already known for all $2 \leq j \leq k-1$.
Setting $u=1$ in \eqref{in_thm_Z-eq-1}--\eqref{in_thm_Z-eq-3}, we obtain, for $1\leq i \leq k-2$,
\vspace{-0.2cm}
\begin{align}\label{in_thm_Z-eq-1-u1}
F_k(x, p,q, v,s,t) = &1 + \sum\limits_{i=0}^{k-2} \Big(f_{k,i}(x, p,q, v,s,t)+ g_{k,i}(x, p,q, v,s,t)\Big),
\\
f_{k,i}(x, p,q, v,s,t)& = pqvx^{i+1} X_i(p,q,1,1,s,t) \left(F_{k-i}(x,p,q,v,1,t) - 1\right) \nonumber  \\
&+
pq^2vx^{i+1} X_i(p,q,1,1,s,1) \left(F_{k-i}(x,p,q,v,1,1) - 1\right)
\left(F_{k}(x,p,q,v,s,t)-1\right)
\nonumber\\
&+
pq
\left(F_{k-i}(x,p,q,v,1,1) - 1\right)
\sum_{j=1}^{i-1}
X_j(p,q,1,1,1,1) x^j   f_{k,i-j}(x, p,q, v,s,t),
\end{align}
\vspace{-0.7cm}
\begin{align} \label{in_thm_Z-eq-3-u1}
g_{k,i}(x, p,q, v,s,t) & =pvtx^{i+1} X_i(p,q,1,1,s,t) \nonumber \\
&+
pqvx^{i+1} X_i(p,q,1,1,s,1)  \left(F_k(x,p,q,v,s,t) - 1\right) \nonumber \\
&+
p^2qvx^{i+1} \sum_{j=1}^{i-1}
  X_j(p,q,1,1,1,1) X_{i-j}(p,q,1,1,s,t)
  \left(F_{k-i+j}(x,p,q,v,1,t) - 1\right) \nonumber\\
  &+
  pq\sum_{j=1}^{i-1} X_{j}(p,q,1,1,1,1) x^j \cdot
\left(
F_{k-i+j}(x,p,q,v,1,1) - 1
\right) \nonumber \\
&~~~~~~~\cdot
\left(
g_{k,i-j}(x, p,q, v,s,t)
-pvtx^{i-j+1}  X_{i-j}(p,q,1,1,s,t)
\right)
\end{align}
with the initial conditions
\begin{align}\label{in_thm_Z-eq-4-u1}
\begin{cases}
f_{k,0}(x, p, q,  v, s, t)& = vstx+qvsx\left(F_k(x, p, q,  v, s, t)-1\right), \vspace{2mm} \\
g_{k,0}(x, p,q,  v,s,t)&= 0, \vspace{2mm}\\
f_{2,0}(x, p, q,  v, s, t) &= \frac{vstx}{1 - qvsx}.
\end{cases}
\end{align}

Step by step, for $i$ from $0$ to $k-2$, it is easy to see that
$f_{k,i}(x, p, q,  v, s, t)$ and
$g_{k,i}(x, p, q,  v, s, t)$
can be expressed by some known factor times $F_k(x, p, q,  v, s, t)$. Substituting these expressions for  $f_{k,i}(x, p, q,  v, s, t)$ and
$g_{k,i}(x, p, q,  v, s, t)$ in
\eqref{in_thm_Z-eq-1-u1},  we obtain the expression for $F_k(x, p, q,  v, s, t)$.
Then, repeating the recursive process in \eqref{in_thm_Z-eq-1}--\eqref{in_thm_Z-eq-3}  with the initial conditions, we finally get the expression of $F_k(x, p, q, u, v, s, t)$.
The proof is completed.
\qed

\begin{cor}
For any $k\geq 2$, we have
	\begin{align}
		F_k(x) = &1 + \sum\limits_{i=0}^{k-2} \Big(f_{k,i}(x)+ g_{k,i}(x)\Big),
		\\
		f_{k,i}(x)=&
		s_{i-1}x^{i+1}  \left(F_{k-i}(x) - 1\right)F_{k}(x) +
		\left(F_{k-i}(x) - 1\right)
		\sum_{j=1}^{i-1}
		s_{j-1} x^j   f_{k,i-j}(x),
	\end{align}
	
	\vspace{-0.7cm}
	
	\begin{align}
		g_{k,i}(x) =&  s_{i-1}x^{i+1} F_k(x) +
		\sum_{j=1}^{i-1} s_{j-1}  x^j \cdot
		\left(
		F_{k-i+j}(x) - 1
		\right)
		g_{k,i-j}(x)
	\end{align}
	for $1\leq i \leq k-2$,
	where $s_n$ is the $n$-th Schr\"oder number, and the initial conditions are
	\begin{align}
		\begin{cases}
			f_{k,0}(x)& = xF_k(x), \vspace{2mm} \\
			g_{k,0}(x)&= 0, \vspace{2mm}\\
			f_{2,0}(x) &= \frac{x}{1 - x}.
		\end{cases}
	\end{align}
\end{cor}

\proof
Setting $p=q=u=v=s=t=1$ in  the formula  for $F_k(x,p,q,u,v,s,t)$ in Theorem~\ref{thm 4}, and then observing that by \cite{West1995},  the number of separable permutations of length $n$ is  $X_n (1,1,1,1,1,1) = s_{n-1}$,
we obtain the desired result.
\qed

\begin{rem}
	\textnormal{
		A limitation of our results in Theorem~\ref{thm 4} is that, for
		$X_n(p,q,u,v,s,t)$, the explicit form of the generating function is known
		only in the case of $X_n(u,v,s,t)$; see ~\cite[Thm~9]{DAM2024CZS}.
		Furthermore, according to ~\cite[Thm~4]{DAM2024CZS}, $S(t,p,q)$ satisfies the cubic algebraic equation:
	}
	\begin{align}\label{eq:pq}
		pq\big(S(t,p,q)\big)^{3} + pq\big(S(t,p,q)\big)^{2} t
		+ S(t,p,q)\big((p+q)t - 1\big) + t = 0,
	\end{align}
	\textnormal{
	Based on the algebraic equation~\eqref{eq:pq}, it is fairly easy to compute
	$X_n(p,q)$ to as many terms as desired. Moreover, by combining these with the
	initial terms of $X_n(u,v,s,t)$, we can also deduce the initial terms of
	$X_n(p,q,u,v,s,t)$.
    } 
\end{rem}

To demonstrate the applications of Theorem \ref{thm 4}, we next provide alternative proofs of Theorems~\ref{F3-thm} and~\ref{dist-general-xyuv-thm}. In Section~\ref{Thm4-cor},  we will also discuss a corollary of Theorem \ref{thm 4}.

\subsection{An alternative proof of Theorem~\ref{F3-thm}}
Letting $k=3$ in Theorem \ref{thm 4} and  using $X_1( p, q, u, v, s, t)=uvst$ and  $F_2(x, p,q, u, v,s,t) =1+\frac{uvstx}{1 - qvsx}$, we obtain
\begin{align*}
F_3(x, p,q,u, v,s,t) &= 1 + \Big(f_{3,0}(x, p,q,u, v,s,t) + g_{3,0}(x, p,q,u, v,s,t)\Big) \nonumber\\
&~~~~~~~~~~~~~+  \Big(f_{3,1}(x, p,q,u, v,s,t) + g_{3,1}(x, p,q,u, v,s,t)\Big), \nonumber
\end{align*}
\vspace{-0.8cm}
\begin{align*}
f_{3,1}(x, p,q,u, v,s,t)
&=pquvx^{2} X_1(p,q,u,1,s,t) \left(F_{2}(x,p,q,1,v,1,t) - 1\right)   \nonumber \\
&~~+
pq^2uvx^{2} X_1(p,q,u,1,s,1) \left(F_{2}(x,p,q,1,v,1,1) - 1\right)  \nonumber \\
&~~~~\cdot
\left(F_{3}(x,p,q,1,v,s,t)-1\right),  \nonumber
\vspace{1mm}\\
&=pquvx^{2} \cdot ust \cdot \frac{vtx}{1 - qvx}   \nonumber\\
&~~~~+
pq^2uvx^{2}  \cdot us \cdot \frac{vx}{1 - qvx}
\left(F_{3}(x,p,q,1,v,s,t)-1\right), \nonumber \\
 \end{align*}
 \vspace{-1.5cm}
 \begin{align*}
 g_{3,1}(x, p,q,u, v,s,t)
 & =puvtx^{2} X_1(p,q,u,1,s,t) \vspace{1mm}   \nonumber \\
&~~~~~+
pquvx^{2} X_1(p,q,u,1,s,1)  \left(F_3(x,p,q,1,v,s,t) - 1\right)  \nonumber
\vspace{1mm}\\
 & =puvtx^{2} \cdot ust +
pquvx^{2} \cdot us  \left(F_3(x,p,q,1,v,s,t) - 1\right),  \nonumber
\end{align*}
\vspace{-1cm}
\begin{align*}
f_{3,0}(x, p, q, u, v, s, t)& = uvstx+quvsx\left(F_3(x, p, q, 1, v, s, t)-1\right), \vspace{1mm}  \nonumber \\
g_{3,0}(x, p,q, u, v,s,t)&= 0.  \nonumber
\end{align*}

Letting $u=1$ in the equations above, we get
\begin{align*}
g_{3,0}(x, p,q, 1, v,s,t) &= 0, \vspace{1mm}\\
f_{3,0}(x, p,q, 1, v,s,t) &= vstx+qvsx\left(F_3(x, p, q,  v, s, t)-1\right), \vspace{1mm}\\
g_{3,1}(x, p,q, 1, v,s,t)  &=  pvst^{2}x^{2}  +
pqv sx^{2}   \left(F_3(x,p,q,v,s,t) - 1\right),
\vspace{1mm}\\
f_{3,1}(x, p,q, 1, v,s,t) & = pqvstx^{2}   \cdot \frac{vtx}{1 - qvx}  +
pq^2vsx^{2}    \cdot \frac{vx}{1 - qvx}
\left(F_{3}(x,p,q,v,s,t)-1\right) \vspace{2mm}\\
& =     \frac{pqv^2st^2x^{3}}{1 - qvx}  +
     \frac{pq^2v^2sx^{3}}{1 - qvx}
\left(F_{3}(x,p,q,v,s,t)-1\right)
\end{align*}
\vspace{-0.8cm}
\begin{align*}
F_3(x, p,q, v,s,t) &= 1 + \Big(f_{3,0}(x, p,q,1, v,s,t) + g_{3,0}(x, p,q,1, v,s,t)\Big)\\
&~~~~~~~~~~~~~+  \Big(f_{3,1}(x, p,q,1, v,s,t) + g_{3,1}(x, p,q,1, v,s,t)\Big)\\
&=1 + vstx+ qvsx\left(F_3(x, p, q, v, s, t)-1\right) \vspace{1mm}
\\
&~~~~
+
\Big(pvst^{2}x^{2}  +
pqv sx^{2}   \left(F_3(x,p,q,v,s,t) - 1\right)
\Big)
\\
&~~~~+
\Big(\frac{pqv^2st^2x^{3}}{1 - qvx} + \frac{pq^2v^2sx^{3}}{1 - qvx}
\left(F_{3}(x,p,q,v,s,t)-1\right)\Big)
\\
&=1 +vstx-qvsx
+pvst^{2}x^{2}  -  pqv sx^{2}
+\frac{pqv^2st^2x^{3}}{1 - qvx} - \frac{pq^2v^2sx^{3}}{1 - qvx}
 \\
&~~~~+\Big(qvsx+pqv sx^{2}+\frac{pq^2v^2sx^{3}}{1 - qvx}
\Big)F_3(x,p,q, v,s,t).
\end{align*}
It follows that
\begin{align*}
F_3(x,p,q, v,s,t)& =\frac{1
+vstx-qvsx
+pvst^{2}x^{2}  -  pqv sx^{2}
+\frac{pqv^2st^2x^{3}}{1 - qvx} - \frac{pq^2v^2sx^{3}}{1 - qvx}}
{1-qvsx-pqv sx^{2}-\frac{pq^2v^2sx^{3}}{1 - qvx}}.
\end{align*}

Therefore, we have
\begin{align*}
F_3(x,p,q,u, v,s,t) &= 1 + \Big(f_{3,0}(x, p,q,u, v,s,t) + g_{3,0}(x, p,q,u, v,s,t)\Big)\\
&~~~~~~~~~~~~~+  \Big(f_{3,1}(x, p,q,u, v,s,t) + g_{3,1}(x, p,q,u, v,s,t)\Big)\\
&=1 +uvstx-quvsx
+pu^2vst^{2}x^{2}  -  pqu^2v sx^{2}
+\frac{pqu^2v^2st^2x^{3}}{1 - qvx} - \frac{pq^2u^2v^2sx^{3}}{1 - qvx}
 \\
&~~~~+\Big(quvsx+pqu^2v sx^{2}+\frac{pq^2u^2v^2sx^{3}}{1 - qvx}
\Big)F_3(x,p,q, v,s,t),
\end{align*}
which is consistent with  Theorem~\ref{F3-thm}.

\subsection{An alternative proof of Theorem~\ref{dist-general-xyuv-thm}}
Taking $k=4$ in Theorem \ref{thm 4}  and using $X_1( p, q, u, v, s, t)=uvst$,
$X_2( p, q, u, v, s, t)=pu^2vst^2+quv^2s^2t$,
and $F_2(x, p,q, u, v,s,t) =1+\frac{uvstx}{1 - qvsx}$,
we obtain the following system of equations.
\begin{align}\label{exam-2-eq-4}
F_4(x, p,q,u, v,s,t) &= 1  + \Big(f_{4,0}(x, p,q,u, v,s,t) + g_{4,0}(x, p,q,u, v,s,t)\Big)\\
&~~~~~+  \Big(f_{4,1}(x, p,q,u, v,s,t) + g_{4,1}(x, p,q,u, v,s,t)\Big)\vspace{1mm}\nonumber \\
&~~~~~+  \Big(f_{4,2}(x, p,q,u, v,s,t) + g_{4,2}(x, p,q,u, v,s,t)\Big),\vspace{1mm} \nonumber
\end{align}

\vspace{-1cm}

\begin{align}
f_{4,0}(x, p, q, u, v, s, t)& = uvstx+quvsx\left(F_4(x, p, q, 1, v, s, t)-1\right), \vspace{1mm} \\
g_{4,0}(x, p,q, u, v,s,t)&= 0,\vspace{1mm}\\
f_{4,1}(x, p,q,u, v,s,t)
&=pquvx^{2} X_1(p,q,u,1,s,t) \left(F_{3}(x,p,q,1,v,1,t) - 1\right)  \\
&~~~~~+
pq^2uvx^{2} X_1(p,q,u,1,s,1) \left(F_{3}(x,p,q,1,v,1,1) - 1\right) \nonumber\\
&~~~~~~~\cdot
\left(F_{4}(x,p,q,1,v,s,t)-1\right),\vspace{1mm} \nonumber\\
&=pquvx^{2} \cdot ust \left(F_{3}(x,p,q,1,v,1,t) - 1\right) \nonumber \\
&~~~~~+
pq^2uvx^{2}  \cdot us \cdot \left(F_{3}(x,p,q,1,v,1,1) - 1\right)
\left(F_{4}(x,p,q,1,v,s,t)-1\right), \vspace{1mm} \nonumber\\
 g_{4,1}(x, p,q,u, v,s,t)
 & =puvtx^{2} X_1(p,q,u,1,s,t) \vspace{1mm}  \\
&~~~~~+
pquvx^{2} X_1(p,q,u,1,s,1)  \left(F_4(x,p,q,1,v,s,t) - 1\right)
\vspace{1mm} \nonumber\\
 & =puvtx^{2} \cdot ust +
pquvx^{2} \cdot us  \left(F_4(x,p,q,1,v,s,t) - 1\right)
\vspace{1mm} \nonumber
\end{align}

\vspace{-1cm}

\begin{align}
f_{4,2}(x, p,q,u, v,s,t)
&=pquvx^{3} X_2(p,q,u,1,s,t) \left(F_{2}(x,p,q,1,v,1,t) - 1\right)  \\
&~~~~~+
pq^2uvx^{3} X_2(p,q,u,1,s,1) \left(F_{2}(x,p,q,1,v,1,1) - 1\right) \nonumber\\
&~~~~~~~\cdot
\left(F_{4}(x,p,q,1,v,s,t)-1\right)\nonumber\\
&~~~~~+pq\left(F_{2}(x,p,q,1,v,1,1) - 1\right)
X_1(p,q,u,1,1,1)xf_{4,1}(x, p,q,u, v,s,t),\vspace{1mm}\nonumber\\
&=pquvx^{3} \cdot (pu^2st^2+qus^2t) \cdot \frac{vtx}{1 - qvx}  \nonumber \\
&~~~~~+
pq^2uvx^{3}  \cdot (pu^2s+qus^2) \cdot
\frac{vx}{1 - qvx} \cdot \left(F_{4}(x,p,q,1,v,s,t) - 1\right) \nonumber\\
&~~~~~+pq \cdot \frac{vx}{1 - qvx} \cdot
uxf_{4,1}(x, p,q,u, v,s,t), \nonumber \vspace{1mm}
\end{align}

\begin{align}\label{g-(4,2)}
 g_{4,2}(x, p,q,u, v,s,t)
 & =puvtx^{3} X_2(p,q,u,1,s,t) \vspace{1mm}  \\
&~~~~~+
pquvx^{3} X_2(p,q,u,1,s,1)  \left(F_4(x,p,q,1,v,s,t) - 1\right)\nonumber\\
&~~~~~+
p^2quvx^{3}
  X_1(p,q,u,1,1,1) X_{1}(x,p,q,u,1,s,t)
  \left(F_{3}(x,p,q,1,v,1,t) - 1\right) \nonumber\\
  &~~~~~+
  pq X_{1}(p,q,u,1,1,1) x\cdot
\left(
F_{3}(x,p,q,1,v,1,1) - 1
\right) \nonumber \\
&~~~~~~~\cdot
\left(
g_{4,1} (x,p,q,u,v,s,t)
-
puvtx^{2}  X_{1}(p,q,u,1,s,t)
\right)
\vspace{1mm}\nonumber\\
 & =
 puvtx^{3} (pu^2st^2+qus^2t) \vspace{1mm}  \nonumber\\
&~~~~~+
pquvx^{3} (pu^2s+qus^2)  \left(F_4(x,p,q,1,v,s,t) - 1\right)\nonumber\\
&~~~~~+
p^2quvx^{3}
  \cdot u  \cdot ust
  \left(F_{3}(x,p,q,1,v,1,t) - 1\right) \nonumber\\
  &~~~~~+
  pq \cdot u  x\cdot
\left(
F_{3}(x,p,q,1,v,1,1) - 1
\right)  \nonumber\\
&~~~~~~~\cdot
\left(
g_{4,1} (x,p,q,u,v,s,t)
-
puvtx^{2}  \cdot ust
\right).\nonumber
\end{align}

Letting $u=1$ in \eqref{exam-2-eq-4}--\eqref{g-(4,2)},
we get
\begin{align*}
F_4(x, p,q, v,s,t) &= 1  + \Big(f_{4,0}(x, p,q,1, v,s,t) + g_{4,0}(x, p,q,1, v,s,t)\Big)\\
&~~~~~+  \Big(f_{4,1}(x, p,q,1, v,s,t) + g_{4,1}(x, p,q,1, v,s,t)\Big)
\vspace{1mm}\\
&~~~~~+  \Big(f_{4,2}(x, p,q,1, v,s,t) + g_{4,2}(x, p,q,1, v,s,t)\Big),
\end{align*}
\vspace{-0.8cm}
\begin{align*}
g_{4,0}(x, p,q, 1, v,s,t)&= 0,
\vspace{1mm}\\
f_{4,0}(x, p, q, 1, v, s, t)& = vstx+qvsx\left(F_4(x, p, q,  v, s, t)-1\right),
\vspace{1mm} \\
 g_{4,1}(x, p,q,1, v,s,t)
 & =pvst^{2}x^{2}  +
pqvsx^{2} \left(F_4(x,p,q,v,s,t) - 1\right)
\vspace{1mm}\\
f_{4,1}(x, p,q,1, v,s,t)
&=pqvstx^{2} \left(F_{3}(x,p,q,v,1,t) - 1\right)  \\
&~~~~~+
pq^2vsx^{2}   \left(F_{3}(x,p,q,v,1,1) - 1\right)
\left(F_{4}(x,p,q,v,s,t)-1\right),
\end{align*}
\begin{align*}
 g_{4,2}(x, p,q,1, v,s,t)
 & =
 pvtx^{3} (pst^2+qs^2t)
 \\
&~~~~~+
pqvx^{3} (ps+qs^2)  \left(F_4(x,p,q,v,s,t) - 1\right)\\
&~~~~~+
p^2qvstx^{3}
  \left(F_{3}(x,p,q,v,1,t) - 1\right) \\
  &~~~~~+
  pqx
\left(
F_{3}(x,p,q,v,1,1) - 1
\right)
\left(
g_{4,1} (x,p,q,1,v,s,t)
-
pvst^{2}x^{2}
\right)
\\
& = p^2vst^3x^{3}+pqvs^2t^2x^3
 \\
&~~~~~+
pqvx^{3} (ps+qs^2)  \left(F_4(x,p,q,v,s,t) - 1\right)\\
&~~~~~+
p^2qvstx^{3}
  \left(F_{3}(x,p,q,v,1,t) - 1\right) \\
  &~~~~~+
  pqx
\left(
F_{3}(x,p,q,v,1,1) - 1
\right)
\cdot pqvsx^{2} \left(F_4(x,p,q,v,s,t) - 1\right),
\end{align*}
\begin{align*}
f_{4,2}(x, p,q,1, v,s,t)
&=pqvx^{3} \cdot (pst^2+qs^2t) \cdot \frac{vtx}{1 - qvx}  \\
&~~~~~+
pq^2vx^{3}  \cdot (ps+qs^2) \cdot
\frac{vx}{1 - qvx} \cdot \left(F_{4}(x,p,q,v,s,t) - 1\right) \\
&~~~~~+pq \cdot \frac{vx}{1 - qvx} \cdot
xf_{4,1}(x, p,q,1, v,s,t)
\vspace{1mm}\\
&= \frac{pqv^2tx^{4}(pst^2+qs^2t)}{1 - qvx}
+ \frac{pq^2v^2x^{4}   (ps+qs^2) }{1 - qvx} \cdot \left(F_{4}(x,p,q,v,s,t) - 1\right) \\
&~~~~~+ \frac{pqvx^2}{1 - qvx}
\cdot
\Big[
pqvstx^{2} \left(F_{3}(x,p,q,v,1,t) - 1\right)  \\
&~~~~~+
pq^2vsx^{2}   \left(F_{3}(x,p,q,v,1,1) - 1\right)
\left(F_{4}(x,p,q,v,s,t)-1\right)
\Big],
\end{align*}
which is consistent with  \eqref{F_4_gf_or_u}.
Hence,
combining with Theorem~\ref{F3-thm}, we obtain the expression of $F_4(x, p,q, v,s,t)$
again, which is the same as \eqref{gf-all-stat-eqn-F4-u}. By verifying that \eqref{exam-2-eq-4}--\eqref{g-(4,2)} is consistent with \eqref{F_4_gf_or},  we complete the proof.

\subsection{A corollary of Theorem~\ref{thm 4}}\label{Thm4-cor}
The study of the joint distribution of the statistics lmax and rmax has been done for separable permutations in \cite{DAM2024CZS}, for permutations avoiding pairs of classical length-3 patterns in \cite{DMTCS2024HS}, and for {\em alternating permutations} in \cite{CarlitzScoville,HanKitZha}. Considering the interest in these statistics in the literature, and the fact that Theorem~\ref{thm 4} is significantly simplified when restricted to them, we state the following corollary as a separate result.

\begin{cor}\label{cor_Y}
For any $k\geq 2$,
the g.f.\ $F_k(x, u, v)$ is given by the following system
\begin{align}\label{thm_3.1-eq}
F_k(x, u, v) &= 1 + \sum\limits_{i=0}^{k-2} \Big(f_{k,i}(x, u, v) + g_{k,i}(x, u, v)\Big), \vspace{1mm}\\
 f_{k,i}(x, u, v) & = \left(F_{k-i}(x,1,v) - 1\right)
\sum_{j=1}^{i}
X_j(u,1) x^j   f_{k,i-j}(x, u, v), \vspace{1mm}\\
 g_{k,i}(x, u, v) & = uvx^{i+1} X_i(u,1)  F_k(x,1,v)
+
  \sum_{j=1}^{i-1} X_{j}(u,1) x^j
\left(
F_{k-i+j}(x,1,v) - 1
\right)
g_{k,i-j}(x,  u, v)
\end{align}
for $1\leq i \leq k-2$,
where the initial conditions are
\begin{align}\label{cor_3.1-eq-2}
f_{k,0}(x, u, v) = uvxF_k(x,1, v), ~~~
g_{k,0}(x, u, v) = 0, ~~~
and~~~f_{2,0}(x, u, v) = \frac{uvx}{1 - vx}.
\end{align}
\end{cor}

\proof
Let $p=q=s=t=1$ in Theorem \ref{thm 4}.
Then, the initial  condition \eqref{in_thm_Z-eq-4} becomes
\begin{align*}
f_{k,0}(x,  u, v)& = uvxF_k(x,  1, v), \vspace{2mm} \\
g_{k,0}(x,  u, v)&= 0, \vspace{2mm}\\
f_{2,0}(x,  u, v) &= \frac{uvx}{1 - vx}.
\end{align*}
Also, equation \eqref{in_thm_Z-eq-1} becomes
\begin{align}\label{in_cor_Z-eq-1}
F_k(x, u, v) = 1 + \sum\limits_{i=0}^{k-2} \Big(f_{k,i}(x, u, v,)+ g_{k,i}(x, u, v)\Big)
\end{align}
 and
\begin{align}\label{in_cor_Z-eq-1-f}
f_{k,i}(x, u, v)
=& uvx^{i+1} X_i(u,1) \left(F_{k-i}(x,1,v) - 1\right)
\nonumber\\
&+
uvx^{i+1} X_i(u,1) \left(F_{k-i}(x,1,v) - 1\right)
\left(F_{k}(x,1,v)-1\right)
\nonumber\\
&+
\left(F_{k-i}(x,1,v) - 1\right)
\sum_{j=1}^{i-1}
X_j(u,1) x^j   f_{k,i-j}(x, u, v)
\nonumber \\
=&
uvx^{i+1} X_i(u,1) \left(F_{k-i}(x,1,v) - 1\right)
F_{k}(x,1,v)
\nonumber\\
&+
\left(F_{k-i}(x,1,v) - 1\right)
\sum_{j=1}^{i-1}
X_j(u,1) x^j   f_{k,i-j}(x, u, v)
\nonumber \\
=&
\left(F_{k-i}(x,1,v) - 1\right)
\sum_{j=1}^{i}
X_j(u,1) x^j   f_{k,i-j}(x, u, v),
\end{align}
for $1\leq i \leq k-2$.

Moreover, equation \eqref{in_thm_Z-eq-3}
becomes
\begin{align} \label{in_cor_Z-eq-3}
g_{k,i}(x,  u, v)
=& uvx^{i+1} X_i(u,1)
+uvx^{i+1} X_i(u,1)  \left(F_k(x,1,v) - 1\right) \nonumber \\
&+
uvx^{i+1} \sum_{j=1}^{i-1}
  X_j(u,1) X_{i-j}(u,1)
  \left(F_{k-i+j}(x,1,v) - 1\right) \nonumber\\
  &+
  \sum_{j=1}^{i-1} X_{j}(u,1) x^j \cdot
\left(
F_{k-i+j}(x,1,v) - 1
\right)
\left(
g_{k,i-j}(x,  u, v)
-uvx^{i-j+1}  X_{i-j}(u,1)
\right)
\nonumber  \\
=& uvx^{i+1} X_i(u,1)  F_k(x,1,v)
\nonumber \\
  &+
  \sum_{j=1}^{i-1} X_{j}(u,1) x^j \cdot
\left(
F_{k-i+j}(x,1,v) - 1
\right)
g_{k,i-j}(x,  u, v)
\end{align}
for $1\leq i \leq k-2$.
The proof is completed.
\qed

The following example, rederiving a particular case of Theorem~\ref{F3-thm}, is an application of Corollary~\ref{cor_Y}.

\noindent
\textbf{Example.}
Taking $k=3$ in \eqref{thm_3.1-eq} and using $X_1(u)=u$, we obtain
\begin{align*}
F_3(x, u, v) &= 1 + \Big(f_{3,0}(x, u, v) + g_{3,0}(x, u, v)\Big)+  \Big(f_{3,1}(x, u, v) + g_{3,1}(x, u, v)\Big),\vspace{1mm}\\
f_{3,1}(x, u, v) & = \big(F_{2}(x, 1, v) - 1\big) \cdot g_{3,1}(x, u, v),\vspace{1mm}\\
 g_{3,1}(x, u, v) & =  X_1(u)x \cdot f_{3,0}(x, u, v)
 =ux  f_{3,0}(x, u, v),
 \vspace{1mm}\\
 f_{3,0}(x, u, v) &= uvxF_3(x, 1, v),
  \vspace{1mm}\\
  g_{3,0}(x, u, v) &= 0.
\end{align*}
Note that $F_2(x, u, v) =1+\frac{uvx}{1-vx}$. Letting $u=1$ in the equations above, we get
\begin{align*}
g_{3,0}(x, 1, v) &= 0, \vspace{1mm}\\
f_{3,0}(x, 1, v) &= vxF_3(x, v), \vspace{1mm}\\
g_{3,1}(x, 1, v)  &=  x  f_{3,0}(x, 1, v)
=vx^2F_3(x, v),\vspace{1mm}\\
f_{3,1}(x, 1, v) & = \big(F_{2}(x, v) - 1\big) \cdot g_{3,1}(x, 1, v)
=\frac{v^2x^3}{1-vx}F_3(x, v).
\end{align*}
Furthermore,
\begin{align*}
F_3(x, v) &= 1 + \Big(f_{3,0}(x,1, v) + g_{3,0}(x,1, v)\Big)+  \Big(f_{3,1}(x,1, v) + g_{3,1}(x, 1,v)\Big)\\
&=1 + vxF_3(x, v)+  \Big(\frac{v^2x^3}{1-vx}F_3(x, v) + vx^2F_3(x, v)\Big)\\
&=1 +\frac{vx-v^2x^2+vx^2}{1-vx}F_3(x, v).
\end{align*}
It follows that
$$F_3(x, v) =\frac{1-vx}{1-2vx-vx^2+v^2x^2}.$$
Therefore, we obtain
\begin{align*}
F_3(x, u, v) &= 1 + \Big(f_{3,0}(x, u, v) + g_{3,0}(x, u, v)\Big)+  \Big(f_{3,1}(x, u, v) + g_{3,1}(x, u, v)\Big)\\
&= 1 +\Big(uvx+u^2vx^2+\frac{u^2v^2x^3}{1-vx}\Big)F_3(x, v) \\
&=1+\frac{uvx+u^2vx^2-uv^2x^2}{1-2vx-vx^2+v^2x^2},
\end{align*}
which is consistent with  Theorem~\ref{F3-thm}.

\section{Concluding remarks}
The main result of this paper, Theorem~\ref{thm 4}, provides a system of functional equations describing the joint distribution of the classical statistics asc, des, lmax, rmax, lmin, and rmin on separable permutations when the POP $P_k$ is avoided for any $k\geq 2$. While our recurrence relations can, in principle, produce the generating functions for the joint distribution for any $k$, we wonder whether there are simpler, more combinatorial approaches to deriving these generating functions. Such approaches should be feasible, at least in cases where additional restrictions are imposed on $S_n(2413,3142,P_k)$. These restrictions would not only reduce the number of permutations under consideration but might also introduce simpler structures, increasing the likelihood of finding connections to other combinatorial objects.

\section*{Acknowledgements.}
The authors are grateful to the anonymous referees for their valuable comments, which helped improve the presentation of the paper. 
The first author is supported by the National Science Foundation of China (No.\ 11801447) and Guangdong Basic and Applied Basic Research Foundation (No.\ 2024A1515011276).

\end{document}